%% file: so_article.tex
\documentclass[review,onefignum,onetabnum]{siamart190516}

\usepackage{subfig}

\input{ex_shared}

\usepackage[commentmarkup=uwave]{changes}

\ifpdf
\hypersetup{
  pdftitle={},
  pdfauthor={}
}
\fi

\begin{document}

\maketitle

\begin{abstract}
  We describe a second-order accurate approach to sparsifying the off-diagonal blocks in the hierarchical approximate factorizations of sparse symmetric positive definite matrices.
  The norm of the error made by the new approach depends quadratically, not linearly, on the error in the low-rank approximation of the given block.
  The analysis of the resulting two-level preconditioner shows that the preconditioner is second-order accurate as well.
  We incorporate the new approach into the recent Sparsified Nested Dissection algorithm
    [\emph{SIAM J. Matrix Anal. Appl., 41 (2020), pp. 715-746}],
  and test it on a wide range of problems.
  The new approach halves the number of Conjugate Gradient iterations needed for convergence, with almost the same factorization complexity, improving the total runtimes of the algorithm.
  Our approach can be incorporated into other rank-structured methods for solving sparse linear systems.
\end{abstract}



\input{intro}

\input{previous}
\input{contributions}

\input{idea}

\input{snd}

\input{experiments}

\input{conclusions}

\input{acknowledgements}

\input{appendix}

\FloatBarrier
\bibliographystyle{siamplain}
\bibliography{references}
\end{document}

%% file: ex_shared.tex
\usepackage{lipsum}
\usepackage{amsfonts}
\usepackage{graphicx}
\usepackage{epstopdf}
\usepackage{algorithmic}
\usepackage{blkarray}
\usepackage{enumitem}
\usepackage{array}
\ifpdf
  \DeclareGraphicsExtensions{.eps,.pdf,.png,.jpg}
\else
  \DeclareGraphicsExtensions{.eps}
\fi

\usepackage[section]{placeins}

\allowdisplaybreaks
\usepackage{tikz,pgfplots}
\usepackage{xcolor}
\usetikzlibrary{patterns}
\usetikzlibrary{quotes,arrows.meta}
\usetikzlibrary{pgfplots.groupplots}
\usepackage{pgfplotstable}

\usepgfplotslibrary{external}
\newcommand{\A}{A}

\newcommand{\tL}{\widetilde{L}}
\newcommand{\hA}{\widehat{A}}
\newcommand{\hS}{\widehat{S}}
\newcommand{\hE}{\widehat{E}}
\newcommand{\tE}{\widetilde{E}}
\newcommand{\hL}{\widehat{L}}
\newcommand{\fo}{first-order scheme}
\newcommand{\fullsec}{full second-order scheme}
\newcommand{\superfine}{superfine second-order scheme}
\newcommand{\Fo}{First-order scheme}
\newcommand{\Fullsec}{Full second-order scheme}
\newcommand{\Superfine}{Superfine second-order scheme}

\definecolor{colorq}{RGB}{44, 162, 95}
\definecolor{colorl}{RGB}{43, 140, 190}
\definecolor{colorc}{RGB}{227, 74, 51}

\definecolor{colory}{RGB}{240, 228, 66}
\definecolor{colorz}{RGB}{0, 114, 178}
\definecolor{colorx}{RGB}{213, 94, 0}
\definecolor{colorb}{RGB}{230, 159, 0}

\pgfplotsset{eps01so/.style={color=colorc, mark=triangle*,mark size=2}}
\pgfplotsset{eps001so/.style={color=colorl, mark=square*,mark size=1.5}}
\pgfplotsset{eps005so/.style={color=colorq, mark=*,mark size=1.5}}
\pgfplotsset{eps01fo/.style={color=colorc, mark=triangle,mark size=2, dash pattern=on 4pt off 1pt on 1pt off 1pt, mark options={solid} }}
\pgfplotsset{eps001fo/.style={color=colorl, mark=square,mark size=1.5, dash pattern=on 4pt off 1pt on 1pt off 1pt, mark options={solid}}}
\pgfplotsset{eps005fo/.style={color=colorq, mark=o,mark size=1.5, dash pattern=on 4pt off 1pt on 1pt off 1pt, mark options={solid}}}

\pgfplotsset{convso1/.style={color=colorc, mark=triangle*,mark size=1}}
\pgfplotsset{convso2/.style={color=colorl, mark=square*,mark size=1}}
\pgfplotsset{convso3/.style={color=colorq, mark=*,mark size=1}}
\pgfplotsset{convso3_nomark/.style={color=colorq, mark=*,mark size=0}}
\pgfplotsset{convfo1/.style={color=colorc, mark=triangle,mark size=1, dash pattern=on 4pt off 1pt on 1pt off 1pt, mark options={solid} }}
\pgfplotsset{convfo2/.style={color=colorl, mark=square,mark size=1, dash pattern=on 4pt off 1pt on 1pt off 1pt, mark options={solid}}}
\pgfplotsset{convfo2_large/.style={color=colorl, mark=square,mark size=2, dash pattern=on 4pt off 1pt on 1pt off 1pt, mark options={solid}}}
\pgfplotsset{convfo3/.style={color=colorq, mark=o,mark size=2, dash pattern=on 4pt off 1pt on 1pt off 1pt, mark options={solid}}}
\pgfplotsset{convso2_round/.style={color=colorl, mark=*,mark size=1}}

\pgfplotsset{convso1g/.style={color=gray, mark=triangle*,mark size=1}}
\pgfplotsset{convfo1g/.style={color=gray, mark=triangle,mark size=2, dash pattern=on 4pt off 1pt on 1pt off 1pt, mark options={solid} }}
\pgfplotsset{convso2b/.style={color=colorb, mark=square*,mark size=1}}
\pgfplotsset{convfo2b/.style={color=colorb, mark=square,mark size=2, dash pattern=on 4pt off 1pt on 1pt off 1pt, mark options={solid}}}
\pgfplotsset{convfos/.style={color=colorc, mark=square,mark size=2, dash pattern=on 4pt off 1pt on 1pt off 1pt}}

\pgfplotsset{convso1g_nomark/.style={color=gray, mark=triangle*,mark size=0}}
\pgfplotsset{convfo1g_nomark/.style={color=gray, mark=triangle,mark size=0, dash pattern=on 4pt off 1pt on 1pt off 1pt, mark options={solid} }}
\pgfplotsset{convso2b_nomark/.style={color=colorb, mark=square*,mark size=0}}
\pgfplotsset{convfo2b_nomark/.style={color=colorb, mark=square,mark size=0, dash pattern=on 4pt off 1pt on 1pt off 1pt, mark options={solid}}}
\pgfplotsset{convfo3_nomark/.style={color=colorq, mark=o,mark size=0, dash pattern=on 4pt off 1pt on 1pt off 1pt, mark options={solid}}}

\newsiamremark{remark}{Remark}
\newsiamremark{hypothesis}{Hypothesis}
\crefname{hypothesis}{Hypothesis}{Hypotheses}
\newsiamthm{claim}{Claim}
\newsiamthm{prop}{Proposition}

\headers{Second Order Accurate Factorization of SPD Matrices}{B. Klockiewicz, L. Cambier, R. Humble, H. Tchelepi, and E. Darve}

\title{Second Order Accurate Hierarchical Approximate Factorization of Sparse SPD Matrices 
}

\author{Bazyli Klockiewicz\thanks{Institute for Computational and Mathematical Engineering, Stanford University
  (\email{bazyli@stanford.edu}, \email{lcambier@stanford.edu}, \email{ryhumble@stanford.edu}).}
\and L\'{e}opold Cambier\footnotemark[2]\and Ryan Humble\footnotemark[2]\and \hspace{5cm} Hamdi Tchelepi\thanks{Department of Energy Resources Engineering, Stanford University
  (\email{tchelepi@stanford.edu}).} \and Eric Darve\thanks{Department of Mechanical Engineering,  Stanford University
  (\email{darve@stanford.edu}).}}

\usepackage{amsopn}


%% file: intro.tex

\section{Introduction}
\label{sec:intro}
Hierarchical approximate factorizations are a recent group of approaches to solving sparse linear systems
\begin{equation*}\label{eq:ax=b}
\begin{matrix}
	Ax = b, & A \in \mathbb{R}^{n \times n}, & b \in \mathbb{R}^n
\end{matrix}
\end{equation*}
such as those arising from discretized partial differential equations (PDEs).
These methods
(also called the \emph{fast hierarchical solvers}) 
include, among others, Hierarchical Interpolative Factorization \cite{ho2016hierarchical,feliu2018recursively}, LoRaSp \cite{pouransari2017fast,yang2019sparse}, or Sparsified Nested Dissection \cite{cambier2020algebraic}.
Significant focus in 
their development
has been on the symmetric positive definite (SPD) case (also  \cite{sushnikova2018compress,chen2018distributed,li2017distributed,chen2019robust}) on which we concentrate in this paper.

Fast hierarchical solvers 
repeatedly sparsify
  selected off-diagonal blocks while performing the block Gaussian elimination of $A$.
Without such sparsifications, Gaussian elimination is not practical for large-scale problems.
To recall why, assume $A$ is a sparse SPD matrix of the form
  \begin{equation*}
  	A =
  	\begin{pmatrix}
  	A_{11} 			& A_{12} \\
  	A_{12}^\top	& A_{22}
  	\end{pmatrix}
  \end{equation*}
  Let $A_{11} = L_1 L_1^\top$ be the Cholesky decomposition. Eliminating the variables corresponding to $A_{11}$ leads to
  \begin{equation*}
  	A = 
  	  \begin{pmatrix}
  	L_{1} 								&  \\
  	A_{12}^\top L_{1}^{-\top}	& I
  	\end{pmatrix}
  	  \begin{pmatrix}
  	I 					& \\
  						& A_{22} - A_{12}^\top A_{11}^{-1} A_{12}
  	\end{pmatrix}
  	  \begin{pmatrix}
  	L_{1}^{\top} 		&L_{1}^{-1} A_{12}   \\
  		& I
  	\end{pmatrix} 
  \end{equation*}
  where the Schur complement $S = A_{22} - A_{12}^\top A_{11}^{-1} A_{12}$ in general contains new fill-in entries and loses the original sparsity of $A_{22}$. 
  If we keep eliminating variables naively, the amount of fill-in will eventually make further computations too expensive, with $\mathcal{O}(n^3)$ complexity of factorizing the matrix.
  The amount of fill-in can be minimized using an ordering of variables, e.g., based on the nested dissection algorithm \cite{george1973nested,lipton1979generalized}.
  Under mild assumptions, nested dissection can reduce the complexity down to $\mathcal{O}(n^{3/2}),$ or $\mathcal{O}(n^2),$ on problems corresponding to, respectively, 2D, and 3D grids.
  However, such complexity is still only practical for relatively small problems.
  
  On the other hand, the fill-in
  arising from elimination typically has the low-rank property.
  Namely, certain off-diagonal 
  blocks in the Schur complement have quickly decaying singular values, and can be sparsified.
To illustrate the sparsification of a single block $S_{12}$ in the Schur complement $S$,
we first scale the leading block 
\[S =   	
\begin{pmatrix}
  	S_{11} 					& S_{12} \\
  	S_{12}^\top			& S_{22}
  	\end{pmatrix} =
   \begin{pmatrix}
  	Z_{1} 			&  \\
  						& I
  	\end{pmatrix} 	
  	\begin{pmatrix}
  	I 							& \hS_{12} \\
  	\hS_{12}^\top			& S_{22}
  	\end{pmatrix}
  	\begin{pmatrix}
  	Z_{1}^\top	&  \\
  						& I
  	\end{pmatrix}
  	= Z S' Z^\top\]
  where $S_{11} = Z_1 Z_1^\top$ is the Cholesky decomposition, and $\hS_{12} = Z_1^{-1} S_{12}$.
Sparsification then involves computing
  a rank-revealing decomposition of $\hS_{12}$, to obtain an orthogonal matrix $Q = \begin{pmatrix} Q_1 & Q_2 \end{pmatrix}$ such that $\| Q_1^T \hS_{12} \| = E $ where $\|E\|$ is small.
 In other words, $Q_2$ approximates the range of $\hS_{12}$.
  Then 
  \begin{align}
  	S' = & \label{eq:sparse-0}
  	  	\begin{pmatrix}
  	I 							& \hS_{12} \\
  	\hS_{12}^\top			& S_{22}
  	\end{pmatrix}  
  	= 
  	\begin{pmatrix}
  	Q_1 & Q_2 \\
  			& 			& I 
  	\end{pmatrix}
   	\begin{pmatrix}
  	I 											&	 										& E \\
  												&	I										& Q_2^\top \hS_{12} \\
  	E^\top		& 	  \hS_{12}^{\top} Q_2			& S_{22}
  	\end{pmatrix}
  	\begin{pmatrix}
  	Q_1^\top 				& \\
  	Q_2^\top				&  \\
  								& I
  	\end{pmatrix} \\ 
  	\approx & \label{eq:sparse-1}
  	\begin{pmatrix}
  	Q & \\
  		& I
  	\end{pmatrix}
  	\underbrace{
   	\begin{pmatrix}
  	I 											&	 										&  \\
  												&	I										& Q_2^\top \hS_{12} \\
  	\											& 	  \hS_{12}^{\top} Q_2	& S_{22}
  	\end{pmatrix}
  	}_{ \begin{matrix} S'' \end{matrix} }
  	\begin{pmatrix}
  	Q^\top 	& \\
  				& I
  	\end{pmatrix}
	 = S' + \mathcal{O}(\| E \|)
  \end{align}
The $E$ and $E^\top$ terms are dropped, and $S''$ is sparser than $S$.
The approaches based on such sparsification scheme include, among others, \cite{ho2016hierarchical,sushnikova2018compress,feliu2018recursively,cambier2020algebraic,chen2018distributed,li2017distributed,yang2019sparse,chen2019robust}. 
These algorithms approximately factorize matrices arising from discretized PDEs, typically in $\mathcal{O}(n)$ or $\mathcal{O}(n \log{n})$ operations, repeatedly sparsifying the matrix during the factorization.
The result is an approximation $A \approx L L^\top$ where $L$ is a product of block-triangular and block-diagonal matrices, so that $L^{-1}L^{-\top} \approx A^{-1}$ can be efficiently applied, typically in $\mathcal{O}(n)$ operations. This operator is then used as a preconditioner in a Krylov subspace method, such as Conjugate Gradient \cite{hestenes1952methods} or GMRES \cite{saad1986gmres}

In this paper, we describe a sparsification approach which makes an $\mathcal{O}(\| E \|^2)$ error in approximating $S'$, while restoring the same amount of sparsity. The general idea is to only drop the $E^\top E$ term when eliminating the leading variables in \cref{eq:sparse-0}:
  \begin{align}
  	S' = & \label{eq:sparse-2}
  	\begin{pmatrix}
  	I 							& \hS_{12} \\
  	\hS_{12}^\top			& S_{22}
  	\end{pmatrix} 
  	= 
  	\begin{pmatrix}
  	Q &  \\
   	E^\top		& I 
  	\end{pmatrix}
   	\begin{pmatrix}
  	I 											&	 										& \\
  												&	I										& Q_2^\top \hS_{12} \\
  			& 	  \hS_{12}^\top Q_2			& S_{22} - E^\top E
  	\end{pmatrix}
  	\begin{pmatrix}
  	Q^\top 				& E \\
  								& I
  	\end{pmatrix} \\ 
  	\approx &
  	\begin{pmatrix}
  	Q & \\
  	E^\top	& I
  	\end{pmatrix}
  	\underbrace{
   	\begin{pmatrix}
  	I 											&	 										&  \\
  												&	I										& Q_2^\top \hS_{12} \\
  	\											& 	  \hS_{12}^\top Q_2	& S_{22}
  	\end{pmatrix}
  	}_{ \begin{matrix} S'' \end{matrix} }
  	\begin{pmatrix}
  	Q^\top 	& E \\
  				& I
  	\end{pmatrix} 
   = S' + \mathcal{O}(\| E \|^2) 
  \end{align}
 The matrix $S''$ is the same as in \cref{eq:sparse-1}.
  	The operator
  	$ \tL^{-\top} \tL^{-1} \approx A^{-1}$ resulting from the approximation, is significantly more accurate than before, but only moderately more expensive to apply (each term $
  	\begin{pmatrix}
  	Q & \\
  	& I
  	\end{pmatrix}$ going into the product $\tL$ is replaced by $	\begin{pmatrix}
  	Q & \\
  	E^\top	& I
  	\end{pmatrix}$ which is sparse, well-conditioned, and whose inverse can be efficiently applied). 

We incorporate the new sparsification approach into the recent Sparsified Nested Dissection algorithm (spaND) \cite{cambier2020algebraic}. The algorithm
 uses the approximation \cref{eq:sparse-1}, which we replace by approximations based on \cref{eq:sparse-2}.
While we focus on the SPD case, in which spaND is guaranteed to complete and be stable without pivoting, the new sparsification approach is applicable also in the general case.

%% file: previous.tex
\subsection{Context}
\label{sec:previous}

Low-rank approximations of the off-diagonal or fill-in matrix blocks is a key ingredient of rank-structured hierarchical methods for solving linear systems that arise from boundary integral equations or discretized PDEs.
These methods include accelerated direct methods based on $\mathcal{H}$-matrices \cite{schmitz2014fast}, hierarchical semi-separable (HSS) matrices \cite{chandrasekaran2007superfast,xia2010fast, xia2010superfast,xia2013efficient}, the hierarchical off-diagonal low-rank framework (HODLR) \cite{aminfar2016fast,kong2011adaptive}, and others \cite{gillman2014direct, amestoy2015improving}. 
At moderate accuracies, these methods can be used as general-purpose preconditioners \cite{grasedyck2008performance,grasedyck2009domain,ghysels2016efficient}.
In particular, \cite{li2012new,xia2010robust,xia2017effective} obtained efficient and robust preconditioners for SPD matrices based on approximation  \cref{eq:sparse-1}. 
A version for general matrices was described in \cite{li2014fast}. 
Hierarchical approximate factorizations \cite{ho2016hierarchical,feliu2018recursively,li2017distributed,pouransari2017fast,sushnikova2018compress,cambier2020algebraic} are recent approaches suitable for sparse systems.
These methods do not require special data-sparse matrix formats that underlie some of the previous approaches. 
 As mentioned above, they explicitly compute an approximate factorization of $A$
 into a product of sparse block-triangular and block-diagonal matrices.

%% file: contributions.tex
\subsection{Contributions}
The contributions of this paper are the following:
\begin{enumerate}
	\item \label{contr-1} We describe the new approach resulting in a quadratic approximation error as in \cref{eq:sparse-2}. 
	We present two variants: a more accurate one (called the \emph{\fullsec}) in which the error term is exactly squared compared to \cref{eq:sparse-1}, and a sparser approach (called \emph{\superfine}) which also has a second-order accuracy.
	\item \label{contr-2} We compute expressions for the condition number of the preconditioned systems for two-level preconditioners resulting from \cref{eq:sparse-1} (which we call the \emph{\fo}), and the new approach. 
	(The two-level preconditioner is obtained by inverting exactly the matrix resulting from sparsifying a single block.) 
	The condition number when using the \fullsec{} depends quadratically, while the condition number when using the \fo{}  depends linearly, on the same term, whose norm is smaller than one. 
	\item \label{contr-3} For the two-level preconditioners, we show that the theoretical bound on the relative error in the preconditioned Conjugate Gradient (PCG) when using the \fullsec{} is exactly squared, compared to the \fo.
	This translates into halving the bound on the iteration count.
	\item \label{contr-4} For right-hand sides that satisfy certain constraints, we also prove that, for the two-level preconditioners used with PCG, the residual at the $k$-th iteration when using the \fullsec{} is exactly the same as the residual at the $2k$-th iteration, when using the \fo.
	As a result the convergence is two times faster at each step of PCG.
	While the constraints are not satisfied exactly in practice, we empirically observe similar behavior in our test cases (see below).
	\item \label{contr-5} We incorporate the new approach into the Sparsified Nested Dissection algorithm (spaND) \cite{cambier2020algebraic}.
	Our benchmarks demonstrate that the new methods involve a \emph{minor cost when computing the preconditioner}.
	\item \label{contr-6} We evaluate the action of the operators $I - MA$ resulting from applying the first- and second-order schemes in spaND (which returns the preconditioning operator $M$), on the eigenvectors of the constant-coefficient Laplace equation.
	We observe that the improvement in the accuracy on most of the spectrum is consistent with the two-level preconditioner analysis. 
	\item \label{contr-7} We perform a study of PCG counts and runtimes when using spaND, as a function of matrix size on high-contrast Laplacians, and run the algorithm on all large SPD matrices from the University of Florida sparse matrix (SuiteSparse) collection \cite{davis2011university}. In all cases, the new approach improves the runtimes of spaND. In particular, consistently among all tested cases, the number of iterations of PCG needed for convergence \emph{is almost exactly halved,} as predicted by the two-level preconditioner analysis.
	\item \label{contr-8} We also observe that, for a given test case and accuracy parameter, the plot of the 2-norm of the residual in PCG as a function of iteration $k$ when using the \fullsec{}, is approximately the same as the one obtained by plotting for $k$ the $2k$-th residual when using the \fo. This is consistent with the theoretical result mentioned in \cref{contr-4}. 
\end{enumerate}


The paper is organized as follows. Our main theoretical results (\cref{contr-1,contr-2,contr-3,contr-4} above) are in
\cref{sec:main}. The description of the spaND algorithm (\cref{contr-5})  is in \cref{sec:snd}. The experimental
results (\cref{contr-6,contr-7,contr-8}) are in \cref{sec:experiments}, followed by conclusions in
\cref{sec:conclusions}.

%% file: idea.tex
\section{First- and second-order approximation schemes using the low-rank property}
We now describe the approximation scheme \cref{eq:sparse-1} and the new approaches resulting in \cref{eq:sparse-2}. 
\label{sec:main}
\subsection{First-order scheme}
\label{sec:main-1}
Consider a sparse SPD matrix of the form
\begin{equation} \label{eq:A}
\A = 
\begin{pmatrix}
			I 				& \A_{12}	&  \\
			\A_{21} 	& \A_{22}	& \A_{23}	 \\
							&	\A_{32}	&	\A_{33}
\end{pmatrix}	
\end{equation}
with a low-rank off-diagonal structure. That is, assume that $\A_{12} = \A_{21}^\top$ has quickly decaying singular values. One can exploit this fact to approximately eliminate a number of variables from the system without introducing any fill-in.
To this end, one computes an orthogonal rank-revealing decomposition of $\A_{12}$ (e.g., the rank-revealing QR or SVD), to obtain a square orthogonal matrix $Q = \begin{pmatrix} Q_f & Q_c \end{pmatrix}$ such that $Q_c$ is a matrix approximating the range of $A_{12}$. In other words, $\|E\| = \mathcal{O}(\varepsilon)$, where $E=Q_f^\top A_{12}$, but $Q_f$ should have as many columns as possible. The first-order scheme is defined by the following approximation, used in \cite{cambier2020algebraic,xia2010robust,xia2017effective,sushnikova2018compress,chen2019robust}
\begin{align}
	\label{eq:Q-basis}
	\A = & \begin{pmatrix}
			Q_f 	& 	Q_c	&\\
			  		& 			& I 		&\\
			  		&			&			& I
			\end{pmatrix}
			\begin{pmatrix}
			I 			&  				 	& E 							&\\
			  			& I 					& Q_c^\top A_{12} 	&\\
			E^\top	& \A_{21} Q_c 	& \A_{22}					& \A_{23}\\
						&						&	\A_{32}					&	\A_{33}	\\								
		   \end{pmatrix}
			\begin{pmatrix}
			Q_f^\top 	& 		& 		\\
			Q_c^\top	& 		&		\\
			 				&	I	&		\\
			 				&		&	I	
			\end{pmatrix} 
			\\ \approx & \label{eq:trailing-1}
			\begin{pmatrix}
			Q_f 	& 	Q_c	&\\
			  		& 			& I 		&\\
			  		&			&			& I
			\end{pmatrix}
			\begin{pmatrix}
			I 			&  										&  \\
			  			& I 										& Q_c^\top A_{12} \\
						& \A_{21} Q_c 	& \A_{22}		& \A_{23}				\\
						&						&	\A_{32}		&	\A_{33}				\\
		    \end{pmatrix}
			\begin{pmatrix}
			Q_f^\top 	& 		& 		\\
			Q_c^\top	& 		&		\\
			 				&	I	&		\\
			 				&		&	I	
			\end{pmatrix}
\end{align}
A number of leading variables (corresponding to $Q_f$) no longer interact with other variables, and are therefore decoupled.
These variables are called the \emph{fine} variables and are denoted by $f$ (the variables corresponding to $Q_c$ are called the \emph{coarse} variables, denoted by $c$). 
The error in approximating $\A$ is given by
\begin{equation}
	\mathcal{E}_1 = 			\begin{pmatrix}
			Q	 	& 		&\\
			  		& I 	&\\
			  		&		& I
			\end{pmatrix}
			\begin{pmatrix}
								&  										& E 	& \\
			  					& 											&     	& \\
			E^\top			& 										 	& 		& \\
								&											&		&
			\end{pmatrix}	
		    \begin{pmatrix}
			Q^\top 	& 		&\\
			 			& I	&\\
			 			&		& I
			\end{pmatrix}
			= \mathcal{O}(\| E \|) = \mathcal{O}(\varepsilon)  
\end{equation}
We note that instead of orthogonal transformations, triangular matrices can also be used (see for example the Hierarchical Interpolative Factorization
\cite{ho2016hierarchical,feliu2018recursively,li2017distributed} or Recursive Skeletonization Factorization \cite{minden2017recursive}).

\subsection{Second-order scheme}
\label{sec:main-2}
The \emph{\fullsec} is obtained by making an error in the Schur complement only, when eliminating the fine variables exactly:
\begin{align}
	\A = 	\label{eq:trailing-0}
			& \begin{pmatrix}
			Q_f		&	Q_c		&		& \\
			E^\top 	&				& I 	& \\
						&				& 		& I
			\end{pmatrix}
			\begin{pmatrix}
			I 			&  										&  								&\\
			  			& I 										& Q_c^\top A_{12}  		&\\
						& \A_{21} Q_c 						& \A_{22} - E^\top E		& \A_{23} \\
						& 											& \A_{32}						& \A_{33}
			\end{pmatrix}	
		    \begin{pmatrix}
			Q_f^\top 		&   E	& \\
			Q_c^\top		&		& \\
			 				& I	& \\
			 				&		& I
			\end{pmatrix} 
			\\ \label{eq:trailing-2}
			\approx & \begin{pmatrix}
			Q_f							&	Q_c		&		& \\
			 \mathbf{E}^\top 	&				& I 	& \\
											&				& 		& I
			\end{pmatrix}
			\begin{pmatrix}
			I 			&  										&  								&\\
			  			& I 										& Q_c^\top A_{12}  		&\\
						& \A_{21} Q_c 						& \A_{22} 					& \A_{23} \\
						& 											& \A_{32}						& \A_{33}
			\end{pmatrix}	
		    \begin{pmatrix}
			Q_f^\top 		& \mathbf{E}	& \\
			Q_c^\top		&		& \\
			 					& I	& \\
			 					&		& I
			\end{pmatrix} 
\end{align}
where we highlighted in bold the new terms appearing in the outer matrices (compared to \cref{eq:trailing-1}).
The error in approximating $A$ this time is
\begin{align}
	\mathcal{E}_2 = &			
	 \begin{pmatrix}
			Q_f							&	Q_c		&		& \\
							 E^\top 	&				& I 	& \\
											&				& 		& I
			\end{pmatrix}
			\begin{pmatrix}
								&  										& \\
			  					& 			-E^\top E 					&     \\
								& 										 	& 
			\end{pmatrix}	
		    \begin{pmatrix}
			Q_f^\top 		& E	& \\
			Q_c^\top	&		& \\
			 				& I	& \\
			 				&		& I
			\end{pmatrix}  \\ =	&		
	\begin{pmatrix}
								&  										& \\
			  					& 			-E^\top E					&  \\
								& 										 	&  
	\end{pmatrix}
	= \mathcal{O}(\| E \|^2) = \mathcal{O}(\varepsilon^2)  
\end{align}
The middle matrix in \cref{eq:trailing-2} is the same as in \cref{eq:trailing-1}. 
It is SPD, if $\A$ is. In fact, its smallest eigenvalue is at least as large as the smallest eigenvalue of the middle matrix in \cref{eq:trailing-0}, i.e., in the exact Cholesky factorization, because a symmetric negative semi-definite term $-E^\top E,$ is dropped. This is called the implicit Schur compensation, and makes the approximation stable \cite{xia2010robust,cambier2020algebraic}.

Notice also that because of the assumed sparsity of $A$, the outer matrices in \cref{eq:trailing-2} involve only a moderate number of additional nonzero entries as compared to the first-order scheme of \cref{eq:trailing-1}.
In a related problem, where $A$ is a dense rank-structured matrix, the \fullsec{} would in general result in a dense factorization (albeit efficiently obtained).


The assumption that the leading block in \cref{eq:A} is the identity matrix, i.e., $\A_{11} = I$, is not limiting. If it is not the case, we scale the leading block
\begin{equation}\label{eq:prescale}
	\A = 
	\begin{pmatrix}
	Z_1	&		&	\\
			& I	&	\\
			&		&	I
	\end{pmatrix}
	\begin{pmatrix}
			I 								& Z_1^{-1}\A_{12}	&  \\
			\A_{21}Z_1^{-\top} 	& \A_{22}					& \A_{23}	 \\
											&	\A_{32}					&	\A_{33}
	\end{pmatrix}
	\begin{pmatrix}
	Z_1^{\top}	&			&	\\
						& I		&	\\
						&			&	I
	\end{pmatrix}	
\end{equation}
where $\A_{11} = Z_1 Z_1^\top$ is the (exact) Cholesky decomposition.
In fact, such prescaling is an essential part of the algorithm, and improves its accuracy and robustness \cite{chen2019robust,cambier2020algebraic,feliu2018recursively}.

\subsubsection{\Superfine}
We can further drop smallest entries of $E$ to obtain a second-order scheme in which the outer matrices are sparser.
The matrix $Q = \begin{pmatrix}
	Q_c &Q_f
\end{pmatrix}$, the set of coarse variables $c$, and the set of decoupled fine variables $f$, remain the same.
However, we further split $Q_f = \begin{pmatrix} Q_{f_1} & Q_{f_2} \end{pmatrix}$ where $Q_{f_1}$ spans the space approximating the left singular vectors of $A_{12}$ whose corresponding singular values are sufficiently small to be immediately neglected (we therefore call $f_1$ the set of \emph{superfine variables}).
More precisely, we can choose $Q_{f_1}$ and $Q_{f_2}$ so that $\| E_1 \| = \mathcal{O}(\varepsilon^2)$, and $\| E_2 \| = \mathcal{O}(\varepsilon)$, where $E_1 = Q_{f_1}^{\top}A_{12}$,  $E_2 = Q_{f_2}^{\top}A_{12}$.
We have 
\begin{align}
	\A = & \begin{pmatrix}
			Q_{f_1} 	& 	Q_{f_2}	& Q_c &\\
			  		& 			& 			& I 		&\\
			  		&			&			& 			& I
			\end{pmatrix}
			\begin{pmatrix}
			I				&						& 										&E_1							&\\
			 				& I 				 	& 										&E_2							&\\
			  				& 						& I 									& Q_c^\top A_{12} 		&\\
			E_1^\top	& E_2^\top		& \A_{21} Q_c 					& \A_{22}						& \A_{23}\\
							&						&										&\A_{32}						&	\A_{33}	\\								
		   \end{pmatrix}
			\begin{pmatrix}
			Q_{f_1}^\top 	& 		& 		\\
			Q_{f_2}^\top	& 		&		\\
			Q_c^\top	& 		&		\\
			 				&	I	&		\\
			 				&		&	I	
			\end{pmatrix} 
			\nonumber \\ \approx & \label{eq:trailing-3}
				 \begin{pmatrix}
			Q_{f_1} 	& 	Q_{f_2}		& Q_c &\\
			  		& 	\mathbf{E}_2		& 			& I 		&\\
			  		&				&			& 			& I
			\end{pmatrix}
			\begin{pmatrix}
			I				&						& 										&									&\\
			 				& I 				 	& 										&							&\\
			  				& 						& I 									& Q_c^\top A_{12} 		&\\
							& 						& \A_{21} Q_c 					& \A_{22}					& \A_{23}\\
							&						&										&\A_{32}						&	\A_{33}	\\								
		   \end{pmatrix}
			\begin{pmatrix}
			Q_{f_1}^\top 	& 		& 		\\
			Q_{f_2}^\top	& 	\mathbf{E}_2	&		\\
			Q_c^\top	& 		&		\\
			 				&	I	&		\\
			 				&		&	I	
			\end{pmatrix}
\end{align}
where we dropped the $E_1$, $E_1^\top$, and $E_2^\top E_2$ terms in the approximation. The middle (trailing) matrix is still the same as in \cref{eq:trailing-1} and \cref{eq:trailing-2}, but the outer matrices are now sparser than in \cref{eq:trailing-2}, while the error is still quadratic: 
\begin{align*}
\mathcal{E}_3 
			& =  \mathcal{O}(\| E_2 \|^2) + \mathcal{O}(\| E_1 \|) = \mathcal{O}(\varepsilon^2)
\end{align*}
The \fo{} can be interpreted as the above scheme with $f_1 = f, f_2 =\emptyset$ while the \fullsec{} is obtained by taking $f_1 = \emptyset, f_2 = f$. The \superfine{} is therefore a ``middle-ground'' scheme that retains, however, the second-order accuracy.

\subsection{Two-level preconditioner analysis}\label{sec:main-analysis}
In practice, the approximations described in \cref{sec:main-1} or \cref{sec:main-2} would be applied recursively in a multilevel algorithm, such as spaND \cite{cambier2020algebraic}, which we describe in \cref{sec:snd}. The spaND algorithm approximately factorizes the matrix to obtain an accurate preconditioner for $A$. The preconditioner can typically be applied in $\mathcal{O}(n)$ operations. To study theoretically the differences between the first- and second-order schemes, we consider a two-level preconditioner, in which the system resulting from the approximation \cref{eq:trailing-1}, or \cref{eq:trailing-2}, or \cref{eq:trailing-3}, is solved exactly. Denote
\begin{equation}
	\hA = 			\begin{pmatrix}
			  			 I 										& Q_c^\top A_{12}  		&\\
						 \A_{21} Q_c 					& \A_{22} 					& \A_{23} \\
																& \A_{32}						& \A_{33}
							\end{pmatrix}	
\end{equation}
The original matrix $A$ can therefore be written as
\begin{equation}\label{eq:Ass}
	\A = ZV
	\begin{pmatrix}
		I 																   & \begin{pmatrix} 0 &E & 0 \end{pmatrix} \\
		\begin{pmatrix} 0 \\ E^\top \\ 0 \\ \end{pmatrix} & \hA 
	\end{pmatrix}
	V^\top Z^\top
\end{equation}
where $Z$ arises from the block-diagonal scaling \cref{eq:prescale}, and $V$ is a sparse orthogonal matrix, such that
 \[ ZV = 	
	\begin{pmatrix}
		Z_1 Q	&		&	\\
				& I	&	\\
				&		&	I
	\end{pmatrix}
\]
We further denote $\widehat{E} = \begin{pmatrix} 0 & E & 0 \end{pmatrix}$.

\subsubsection{\Fo}
In the first-order scheme we drop the terms $\widehat{E}$ and $\widehat{E}^{\top}$ in \cref{eq:Ass}. The two-level preconditioner has the form
$
	M_1 =L_1L_1^\top
$
where
\begin{equation}\label{eq:l1}
	L_1 = ZV
\begin{pmatrix}
	   I  	&	\\
	   		&	\hL
\end{pmatrix}
\end{equation}
with $\hA=\hL \hL^\top$ being the (exact) Cholesky decomposition of $\hA$.
\begin{prop} \label{prop:first}
(See \cite{xia2010robust}, Proposition 2.1)
Let $A$ be an SPD matrix, and let $ M_1 = L_1 L_1^\top$ be the preconditioner defined by the first-order scheme, with $L_1$ as in \cref{eq:l1}. Then
\[ 
L_1^{-1}AL_1^{-\top}  =
\begin{pmatrix}
	  I								& \hE \hL^{-\top} \\
	  \hL^{-1}\hE^{\top}	 	& I	
\end{pmatrix} \]
where $\|\hL^{-1}\hE^\top\|_2 < 1.$ In particular, the 2-norm condition number of the preconditioned system is given by
\begin{equation}\label{eq:cond-1}
	 \kappa(L_1^{-1}AL_1^{-\top}) = \frac{1 + \| \hL^{-1} \hE^\top \|_2}{1 - \| \hL^{-1} \hE^\top \|_2}
\end{equation}
\end{prop}
\begin{proof}
We have 
\begin{equation}
	L_1^{-1}AL_1^{-\top} = 
\begin{pmatrix}
	   I  	&	\\
	   		&	\hL^{-1}
\end{pmatrix}
\begin{pmatrix}
	  I								& \hE \\
	  \hE^{\top}				& \hA	
\end{pmatrix}
\begin{pmatrix}
	   I  	&	\\
	   		&	\hL^{-\top}
\end{pmatrix}
=
\begin{pmatrix}
	  I								& \hE \hL^{-\top} \\
	  \hL^{-1}\hE^{\top}	 	& I	
\end{pmatrix}
\end{equation}
 Notice that $\| \hL^{-1} \hE^{\top} \|_2 < 1$ because the Schur complement $I - \hL^{-1}\hE^\top \hE \hL^{-\top}$ is SPD. The condition number of a matrix
 $\begin{pmatrix}
 	I 			& C \\
 	C^\top	& I
 \end{pmatrix}$ with $\| C \|_2 < 1$ equals $\frac{1 + \|C\|_2}{1 - \|C\|_2}.$
\end{proof}
\subsubsection{\Fullsec}
We have
\begin{equation}
A =
ZV
\begin{pmatrix}
  I  	     & \\
  \hE^\top & I
\end{pmatrix}
\begin{pmatrix}
	   I  	&	\\
	   		& \hA - \hE^\top \hE
\end{pmatrix}
\begin{pmatrix}
  I  	     & \hE \\
   			 & I
\end{pmatrix}
V^\top Z^\top
\end{equation}
In the \fullsec, we only drop the term $\hE^\top \hE$ above.
 The two-level preconditioner therefore has the form
 $
 	M_2 = L_2L_2^\top
$
 where
 \begin{equation}\label{eq:l2}
 L_2 = ZV
\begin{pmatrix}
  I  	     	& \\
  \hE^\top & I
\end{pmatrix}
\begin{pmatrix}
	   I  	&	\\
	   		& \hL
\end{pmatrix}
 \end{equation}

\begin{prop}\label{prop:second}
	Let $A$ be an SPD matrix, and let $M_2=L_2 L_2^\top$ be the preconditioner defined by the \fullsec, with $L_2$ as in \cref{eq:l2}. Then
	\[ L_2^{-1}AL_2^{-T} =  
	\begin{pmatrix}
				 I	&  \\
					& I - \hL^{-1} \hE^\top \hE \hL^{-\top}
	\end{pmatrix} 
	\] 
	where $ \| \hL^{-1} \hE^\top \|_2 < 1$. In particular, the 2-norm condition number of the preconditioned system is given by
	\begin{equation}\label{eq:cond-2}
		\kappa(L_2^{-1} A L_2^{-T}) = \frac{1}{1 - \| \hL^{-1}\hE^\top \|^{{\mathbf{2}}}_2}
	\end{equation}
\end{prop}
\begin{proof}
	
We compute
\begin{align*}
	& L_2^{-1} A L_2^{-\top} =
\begin{pmatrix}
	   I  	&	\\
	   		& \hL^{-1}
\end{pmatrix}
\begin{pmatrix}
  I  	     		& \\
  -\hE^\top 	& I
\end{pmatrix}	
\begin{pmatrix}
	  I								& \hE \\
	  \hE^{\top}				& \hA	
\end{pmatrix} 
\begin{pmatrix}
  I  	     		& -\hE \\
  				 	& I
\end{pmatrix}
\begin{pmatrix}
	   I  	&	\\
	   		& \hL^{-\top}
\end{pmatrix} \\
	& = 
\begin{pmatrix}
	   I  	&	\\
	   		& \hL^{-1}
\end{pmatrix}
	\begin{pmatrix}
	   I  	&	\\
	   		& \hA - \hE^\top \hE
\end{pmatrix}
\begin{pmatrix}
	   I  	&	\\
	   		& \hL^{-\top}
\end{pmatrix} =   I -  	\begin{pmatrix}
	     	&	\\
	   		& \hL^{-1} \hE^\top \hE \hL^{-\top}
\end{pmatrix}
\end{align*}
 Since this matrix is also SPD, we confirm that $ \| \hL^{-1} \hE^\top \|_2 < 1$. To obtain the expression for the condition number, notice that the smallest eigenvalue of $L_2^{-1}AL_2^{-\top}$ equals $1 - \| \hL^{-1} \hE^\top \|_2^2,$ and the largest one equals 1.
\end{proof}

 Comparing \Cref{prop:first} and \Cref{prop:second} we can see that the second-order scheme is strictly more accurate than the first-order scheme in terms of the preconditioner accuracy. The new terms
 $\begin{pmatrix}
  I  	     		& \\
  - \hE^\top 			 & I
\end{pmatrix}$ and 
$\begin{pmatrix}
  I  	     		& - \hE \\
  		 			& I
\end{pmatrix}$
can be efficiently applied to a vector. The Taylor series expansions $ \frac{1+x}{1-x} = 1 + 2x + O(x^2) $ and $\frac{1}{1-x^2} = 1 + x^2 + O(x^4) $  at $x=0$, justify the term ``second-order''. Notice also that $ \| \hL^{-1} \hE^\top \|_2 < 1$ would be true for any choice of orthogonal $Q = \begin{pmatrix} Q_c & Q_f \end{pmatrix}$ (even if, for example, we were to maximize $\| \hE \|_2$ instead of minimizing it).

The scaling matrix $Z$ can prescale more blocks than just $A_{11}.$
It can be a matrix prescaling all diagonal blocks ahead of time (this is called the double-sided scaling \cite{xia2011some,xia2017effective}).
In that case, $Z$ is a block Jacobi preconditioner which preconditions the matrix $ \hA.$
From \Cref{prop:first} and \Cref{prop:second}, choosing $Z$ in this way should improve the preconditioner quality in both the first- and second-order schemes.
For the \fo, this was demonstrated in \cite{feliu2018recursively,cambier2020algebraic,xia2017effective}.

\subsubsection{\Superfine}
The preconditioner resulting from the \superfine{} is given by
$
	M_3 = L_3 L_3^\top
$
with
\begin{equation}\label{eq:l3}
	L_3 = ZV 
	\begin{pmatrix}
	I 				& 		&  \\
					& I	&  \\
					&	\hE_2^\top		& I \\
	\end{pmatrix}
	\begin{pmatrix}
	   I  	&	\\
	   		& \hL
\end{pmatrix}
\end{equation}
where $\hE_2 = \begin{pmatrix} 0 & E_2 & 0 \end{pmatrix}.$ Denote also $\hE_1 = \begin{pmatrix} 0 & E_1 & 0 \end{pmatrix}.$ Recall that $\| \hE_2 \|_2 = \mathcal{O}(\varepsilon)$ and $\| \hE_1 \|_2 = \mathcal{O}(\varepsilon^2).$ With a similar computation as before, we obtain
\begin{prop}\label{prop:second-sf}
	Let $A$ be an SPD matrix, and let $M_3=L_3 L_3^\top$ be the preconditioner defined by the \superfine, with $L_3$ as in \cref{eq:l3}. Then
\[
	L_3^{-1} A L_3^{-\top} = I + 
\begin{pmatrix}
	 								& 							& \hE_1 \hL^{-\top}  \\
									& 							&  \\
	\hL^{-1} \hE_1^\top	&							&  \\
\end{pmatrix}
-
\begin{pmatrix}
	 					& 							&  \\
						& 							&  \\
						&							& \hL^{-1} \hE_2^\top \hE_2 \hL^{-\top} \\
\end{pmatrix}
\]
\end{prop}
\begin{proof} We have:
\begin{align*}
	 L_3^{-1} A L_3^{-\top} & =
\begin{pmatrix}
	   I  	&	\\
	   		& \hL^{-1}
\end{pmatrix}
		\begin{pmatrix}
	I 					& 							& \hE_1  \\
						& I						& \\
	\hE_1^\top	&		& \hA  -\hE_2^\top \hE_2 \\
	\end{pmatrix}
\begin{pmatrix}
	   I  	&	\\
	   		& \hL^{-\top}
\end{pmatrix} \\
& =
\begin{pmatrix}
	I 								& 							& \hE_1 \hL^{-\top} \\
									& I						&  \\
	\hL^{-1} \hE_1^\top	&							& I- \hL^{-1} \hE_2^\top \hE_2 \hL^{-\top} \\
\end{pmatrix} 
\end{align*}
	
\end{proof}
Clearly, \Cref{prop:second-sf} contains \Cref{prop:first} and \Cref{prop:second} as special cases. 
\subsubsection{The bound on iteration count is halved}
The error made by Conjugate Gradient in solving $Ax = b$ is bounded by \cite{trefethen1997numerical,golub2012matrix}
\begin{equation}\label{eq:cg-convergence}
	\| x - x_k \|_A  \leq 2 \| x - x_0 \|_A \left(\frac{\sqrt{\kappa(A)} - 1}{\sqrt{\kappa(A)} + 1}\right)^k 
\end{equation}
where $x$ is the exact solution, $x_k$ is the approximate solution at the $k$-th iteration, and $\kappa(A)$ is the 2-norm condition number of $A$.
The term $R =  ({\sqrt{\kappa(A)} - 1})/({\sqrt{\kappa(A)} + 1})$ is thus a bound on convergence rate.
Denote the preconditioned systems after applying the first- and the full second-order schemes by, respectively, $A_1 = L_1^{-1} A L_1^{-\top}$ and $A_2 = L_2^{-1} A L_2^{-\top}.$
From \cref{eq:cond-1} and \cref{eq:cond-2}, denoting $\tilde{\varepsilon} = \| \hL^{-1} \hE^\top \|_2$ the corresponding bounds on convergence rates in the respective norms $\| \ . \ \|_{A_1}$ and $\| \ . \ \|_{A_2},$ are given by
\begin{displaymath}
	R_1 = \frac{\sqrt{\frac{1 + \tilde{\varepsilon}}{1 - \tilde{\varepsilon}}} - 1 }{\sqrt{\frac{1 + \tilde{\varepsilon}}{1 - \tilde{\varepsilon}}} + 1 }
	\ \text{ and }
	R_2 = \frac{\sqrt{\frac{1}{1 - \tilde{\varepsilon}^2}} - 1 }{\sqrt{\frac{1}{1 - \tilde{\varepsilon}^2}} + 1 }
\end{displaymath}
Now, notice that since $ 0 \leq \tilde{\varepsilon} < 1$ we have
\begin{displaymath}
	R_1^2 = \left(  \frac{\sqrt{1 +  \tilde{\varepsilon}} - \sqrt{1 -\tilde{\varepsilon}}}{\sqrt{1 +  \tilde{\varepsilon}} + \sqrt{1 -\tilde{\varepsilon}}} \right)^2 = 
	\frac{1 - \sqrt{1 - \tilde{\varepsilon}^2}}{1 + \sqrt{1 - \tilde{\varepsilon}^2}} =
	\frac{\frac{1}{ \sqrt{1 - \tilde{\varepsilon}^2}} - 1}{\frac{1}{\sqrt{1 - \tilde{\varepsilon}^2}} + 1}
	= R_2
\end{displaymath}
so the \fullsec{} bound on convergence rate is an exact square of the \fo{} bound on convergence rate.
Let $k_1$ and $k_2$ denote the numbers of iterations needed for convergence to a specified tolerance in the respective norms, for the first- and the \fullsec, that is 
$
	\frac{\| x - x_{k_1} \|_{A_1}}{ \| x - x_0 \|_{A_1}}  = \frac{\| x - x_{k_2} \|_{A_2}}{\| x - x_0 \|_{A_2}}  = \eta
$
where $0 < \eta < 1$ is a target small number. Assuming \cref{eq:cg-convergence} is tight, we obtain
\begin{displaymath}
	R_1^{k_1} = R_2^{k_2} = R_1^{2 k_2} \implies k_2 = \frac{k_1}{2}
\end{displaymath}
The bound on the number of iterations needed for convergence is therefore exactly halved (notice also that the norm $\| \ . \ \|_{A_2}$ is closer to $\| \ . \ \|_2$ than $\| \ . \ \|_{A_1}$).
\subsubsection{Convergence in a special case}
\label{sec:thm}
For right-hand side vectors satisfying a certain constraint, the convergence of Conjugate Gradient is exactly two times faster at each step in the algorithm.
Namely, the following result holds, whose proof we include in \cref{sec:proof}.
\begin{theorem}
\label{thm:residuals}
Consider the Conjugate Gradient algorithm applied to systems such as those obtained when using the two-level preconditioner with the \fo{} and the \fullsec{}, respectively:
\begin{equation}\label{eq:cg-prec-1}
A_1 x = 
\begin{pmatrix}
I & C \\
C^T & I
\end{pmatrix}
\begin{pmatrix}
	x_1 \\ x_2
\end{pmatrix}
=
\begin{pmatrix}
	b_1 \\ b_2
\end{pmatrix}
= b
\end{equation}
\begin{equation}\label{eq:cg-prec-2}
	A_2 x =
\begin{pmatrix}
I & 0 \\
0 & I - C^T C
\end{pmatrix}
\begin{pmatrix}
	x_1 \\ x_2
\end{pmatrix}
= 
\begin{pmatrix}
	b_1 \\ b_2
\end{pmatrix}
= b
\end{equation}
where $\|C\|_2 < 1$ and we additionally assume that
\begin{equation}
	C^\top b_1 = 0
\end{equation}
Let $r^{(1)}_k = b - A_1 x^{(1)}_k$, and $r^{(2)}_k= b - A_2 x^{(2)}_k$ denote the residuals at the $k$-th step of the algorithm run with a zero initial guess $x_0 = 0$, applied to, respectively, \cref{eq:cg-prec-1}, and \cref{eq:cg-prec-2}. Then:
\[ r^{(1)}_{2k} = r^{(2)}_{k}  \]
i.e., the residual produced by Conjugate Gradient at the $2k$-th iteration, when applied to \cref{eq:cg-prec-1}, equals the residual produced at the $k$-th iteration, when applied to \cref{eq:cg-prec-2}.
\end{theorem}
Below, we show example plots obtained from small problems satisfying the constraint $C^\top b_1 = 0$ exactly, illustrating the theorem.
In a model two-level preconditioner, while $C^\top = \hL^{-1} \hE^\top$ may not satisfy the constraint exactly for a given right-hand side, it has small singular values quickly decaying to zero.
See \cref{sec:experiments} for empirical plots obtained with the spaND algorithm.

\pgfplotstableread[col sep=comma]{logs/first-order-example-convergence.csv}\dataex
\pgfplotstableread[col sep=comma]{logs/second-order-example-convergence.csv}\dataexx
\begin{figure}[htbp]
    \centering
    \begin{tikzpicture}
        \begin{groupplot}[
            group style={
                group name=cg_example,
                group size=1 by 1,
                xlabels at=edge bottom,
                xticklabels at=edge bottom,
                vertical sep=0.2cm,
                horizontal sep=0.2cm,
                ylabels at=edge left,
                yticklabels at=edge left,
            },
            ymode=log,ymin=7*10^(-7),ymax=50.0,
            ytick={10^2,10^0,10^(-2),10^(-4),10^(-6)},
            yticklabels={$10^2$,$10^0$,$10^{-2}$,$10^{-4}$,$10^{-6}$},
            xmin=1,xmax=30,
            xtick={10,20,30},
            xticklabels={10,20,30},
            extra y ticks={16.38538372,2.718902135,0.610112573,0.151672647,0.042175844,0.013156101,0.003805835,0.001084392,0.000311419,9.52E-05,2.91E-05,7.93E-06,2.89E-06,9.1E-07},
            extra y tick labels={},
            extra tick style={
            grid=major,
        },
        ]
        \nextgroupplot[width=10cm,height=5cm,,xlabel={CG iteration, $k$},
        	ylabel style={align=center},
        	ylabel={Residual norm},
        	legend columns = 1,
        	legend entries={First-Order $\| r_k^{(1)} \|_2$, Second-Order $\| r_k^{(2)} \|_2$, First-Order $\|r^{(1)}_{2k}\|_2$ },
        	legend style={font=\footnotesize}]
             \addplot[convfo2] table[x=it, y=res1]{\dataex};
             \addplot[convso2_round] table[x=it, y=res2]{\dataexx};
			 \addplot[convfos] table[x=it, y=res1s]{\dataexx};
        \end{groupplot}
    \end{tikzpicture}
    \caption{Illustration of \cref{thm:residuals} }
    \label{fig:cg-example}
\end{figure}


\subsection{Contrast with related work}
Our work is closely related to that of Xia \cite{xia2010fast,xia2017effective,xia2013efficient}, Li \cite{li2012new,li2014fast}, and Xin \cite{xin2020effectiveness} (which concentrate however on the general, not sparse, case).
In particular, our approach includes the double-sided scaling proposed in those works, as well as the implicit Schur compensation.
The most significant difference is that we include an additional term $E$ in \cref{eq:trailing-2} (or $E_2$ in \cref{eq:trailing-3}) after performing an explicit change of basis, which yields the quadratic approximation error.
More precisely, the error in our approach is the $\mathcal{O}(\varepsilon^2)$ error resulting from dropping the Schur complement update, and possibly an $\mathcal{O}(\varepsilon^2)$ error in the elimination matrices.
Thus, the overall approximation error is also $\mathcal{O}(\varepsilon^2)$ while the other approaches have an $\mathcal{O}(\varepsilon)$ overall approximation error.
Lastly, our methods do not require any special low-rank matrix formats, such as the HSS format used by some of the related approaches.

Most recently, Xia \cite{xia2020robust} used similar ideas, to improve the SIF algorithm for general (typically dense) SPD matrices \cite{xin2020effectiveness}.
Below, we use our methods to improve spaND \cite{cambier2020algebraic}, which is an algorithm for sparse matrices, related to, but different from SIF.

%% file: snd.tex
\section{Hierarchical multilevel algorithm}
\label{sec:snd}

Sparsified Nested Dissection (spaND) \cite{cambier2020algebraic} is a hierarchical multilevel algorithm which repeatedly applies the first-order approximation scheme from \cref{sec:main-1}.
The first-order scheme can be easily replaced by the second-order schemes of \cref{sec:main-2}.
The spaND algorithm is guaranteed to complete on any SPD matrix.
The result is an approximate factorization of $A$ which---under typical assumptions on the ranks of the fill-in blocks---can be computed in $\mathcal{O}(n \log(n))$ operations on a sparse matrix such as those arising from discretized PDEs.
The resulting preconditioner can then be applied in $\mathcal{O}(n)$ operations.

\begin{figure}
    \centering
    \subfloat[][Partitioning defined on the graph of $A$]{\includegraphics[width=0.3\textwidth]{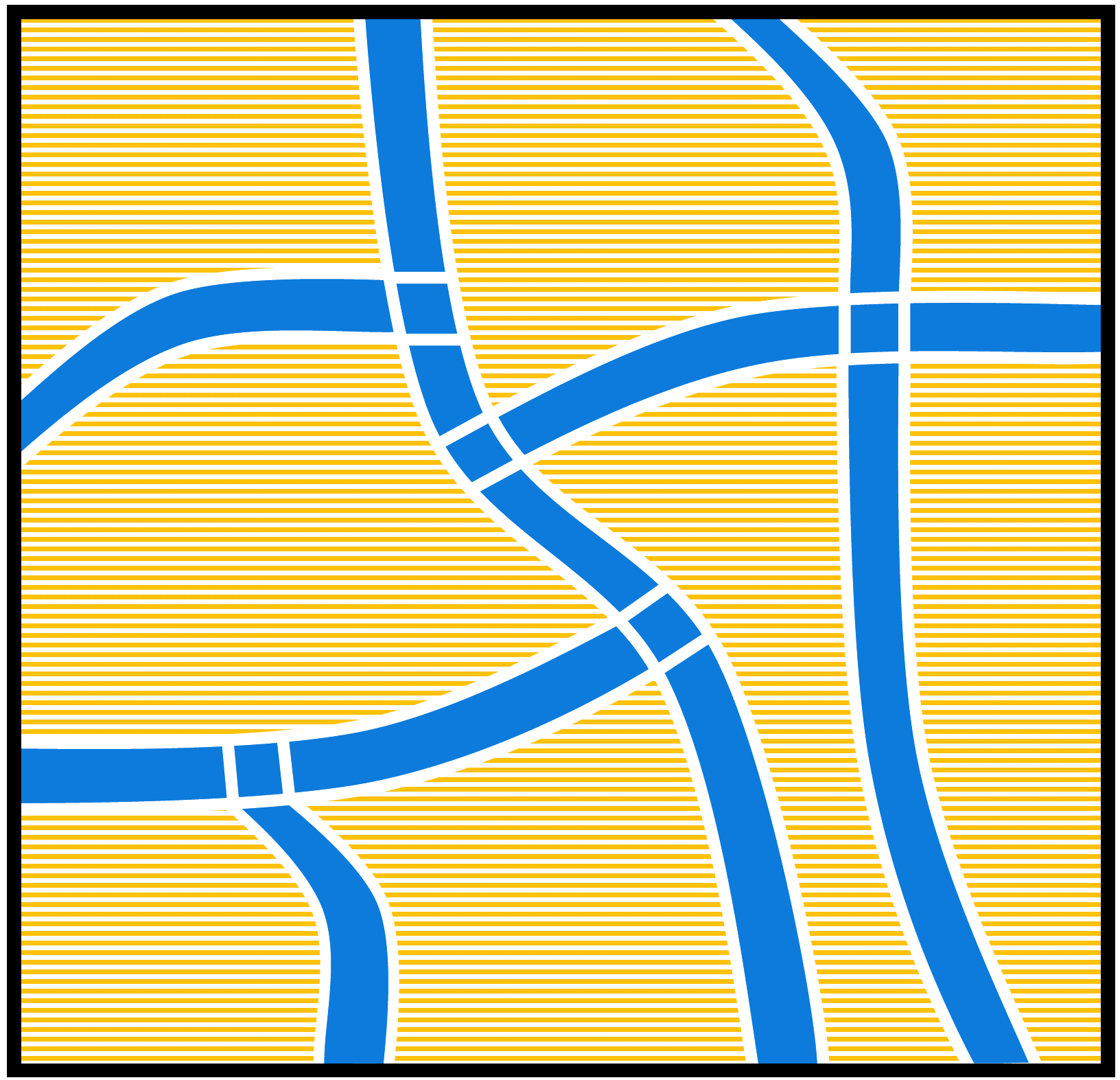}} \quad
    \subfloat[][$l=3$, after eliminating interiors]{\includegraphics[width=0.3\textwidth]{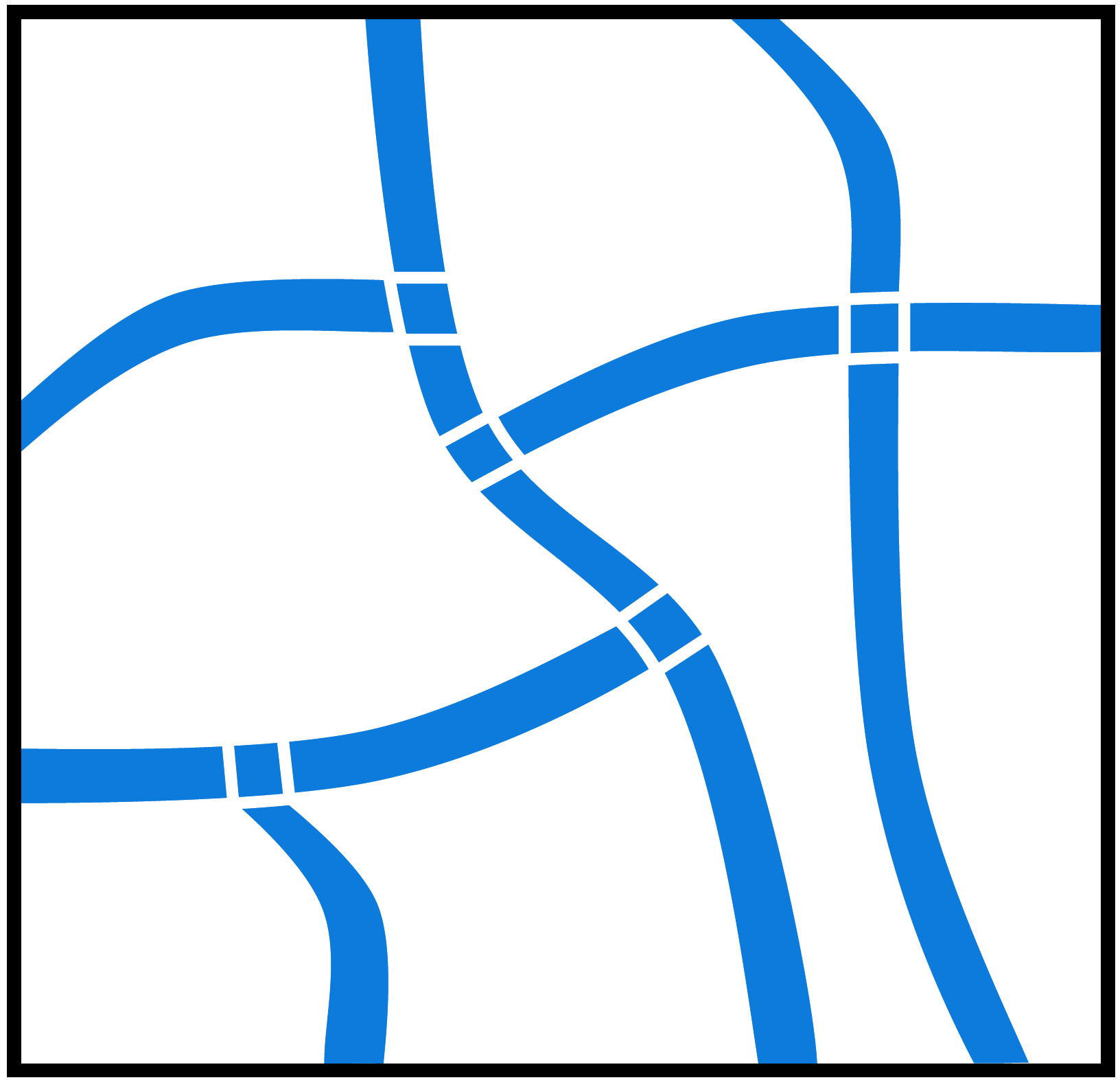}} \quad
    \subfloat[][$l=3$, after sparsifying interfaces]{\includegraphics[width=0.3\textwidth]{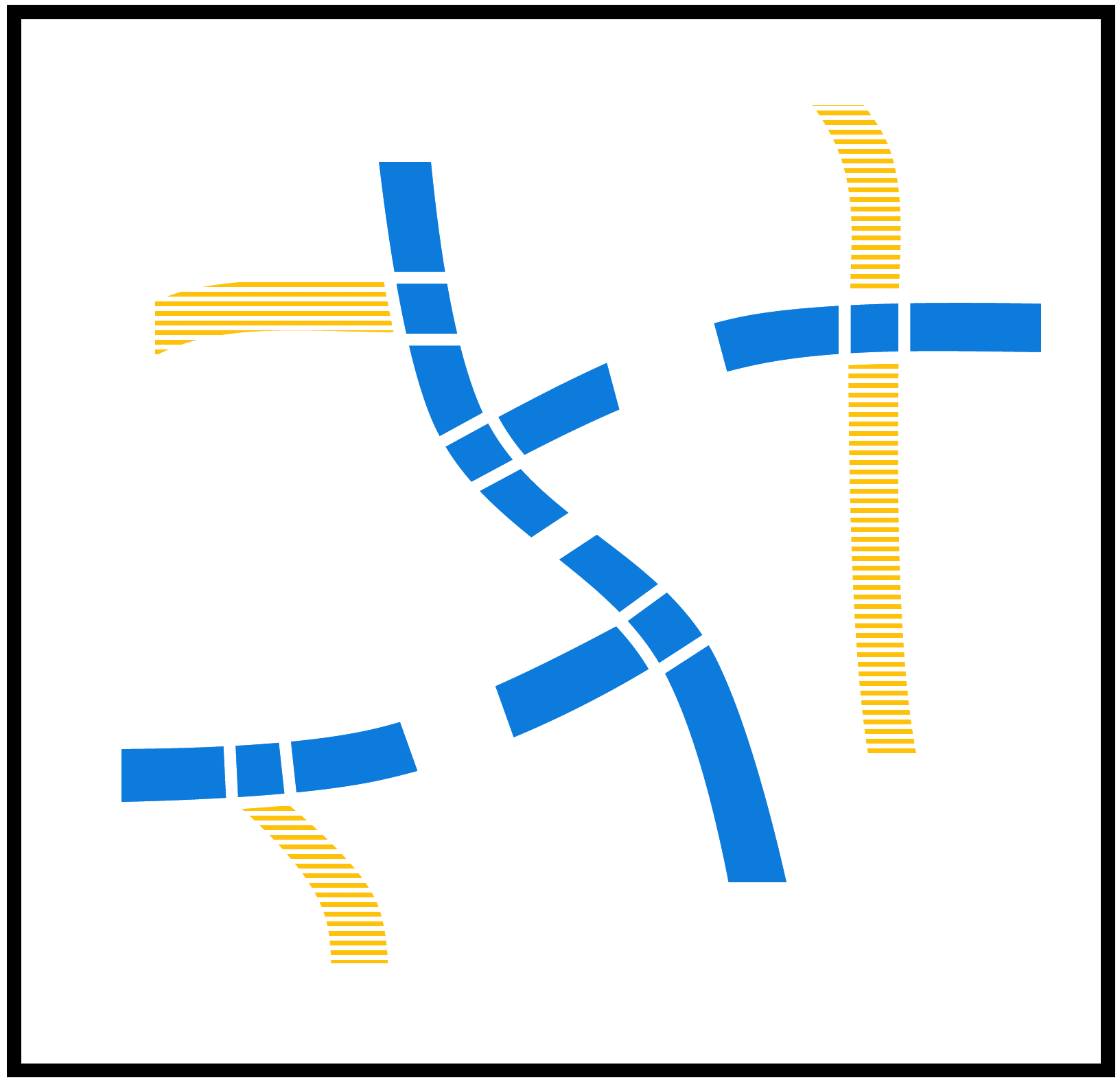}} \\
    \subfloat[][$l=2$, after eliminating interiors]{\includegraphics[width=0.3\textwidth]{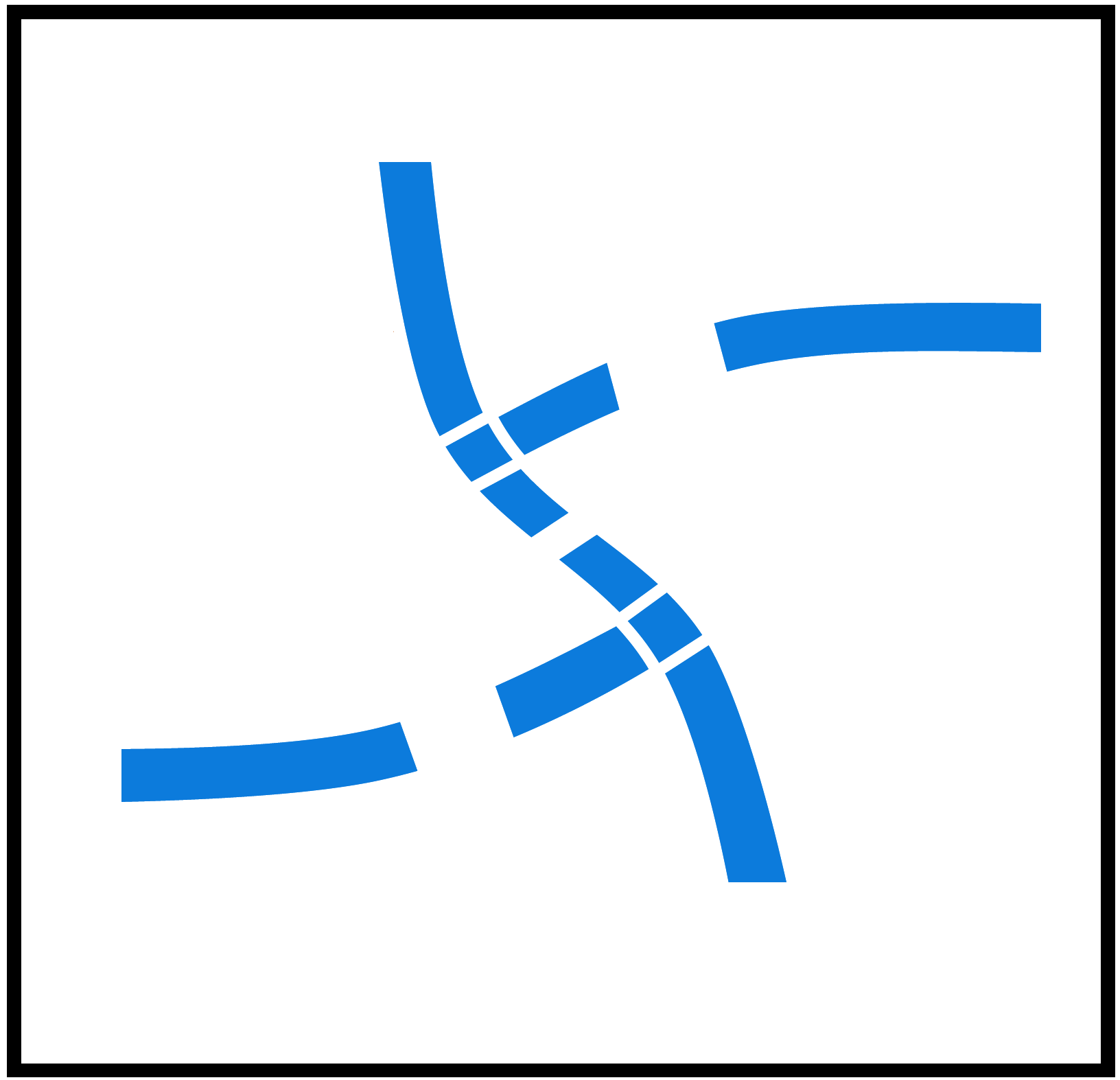}} \quad
    \subfloat[][$l=2$, after sparsifying interfaces]{\includegraphics[width=0.3\textwidth]{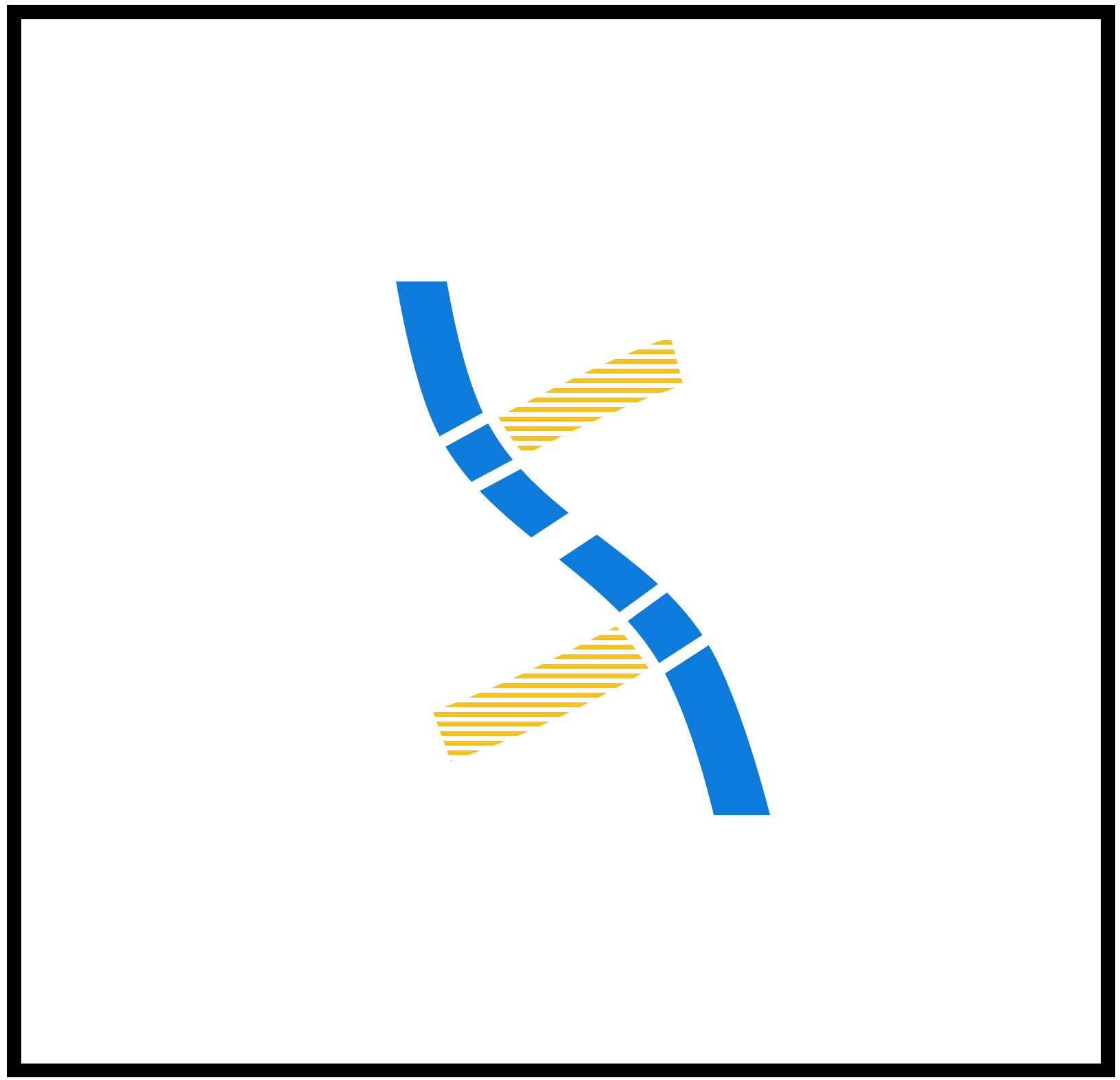}} \quad
    \subfloat[][$l=1$, after eliminating interiors]{\includegraphics[width=0.3\textwidth]{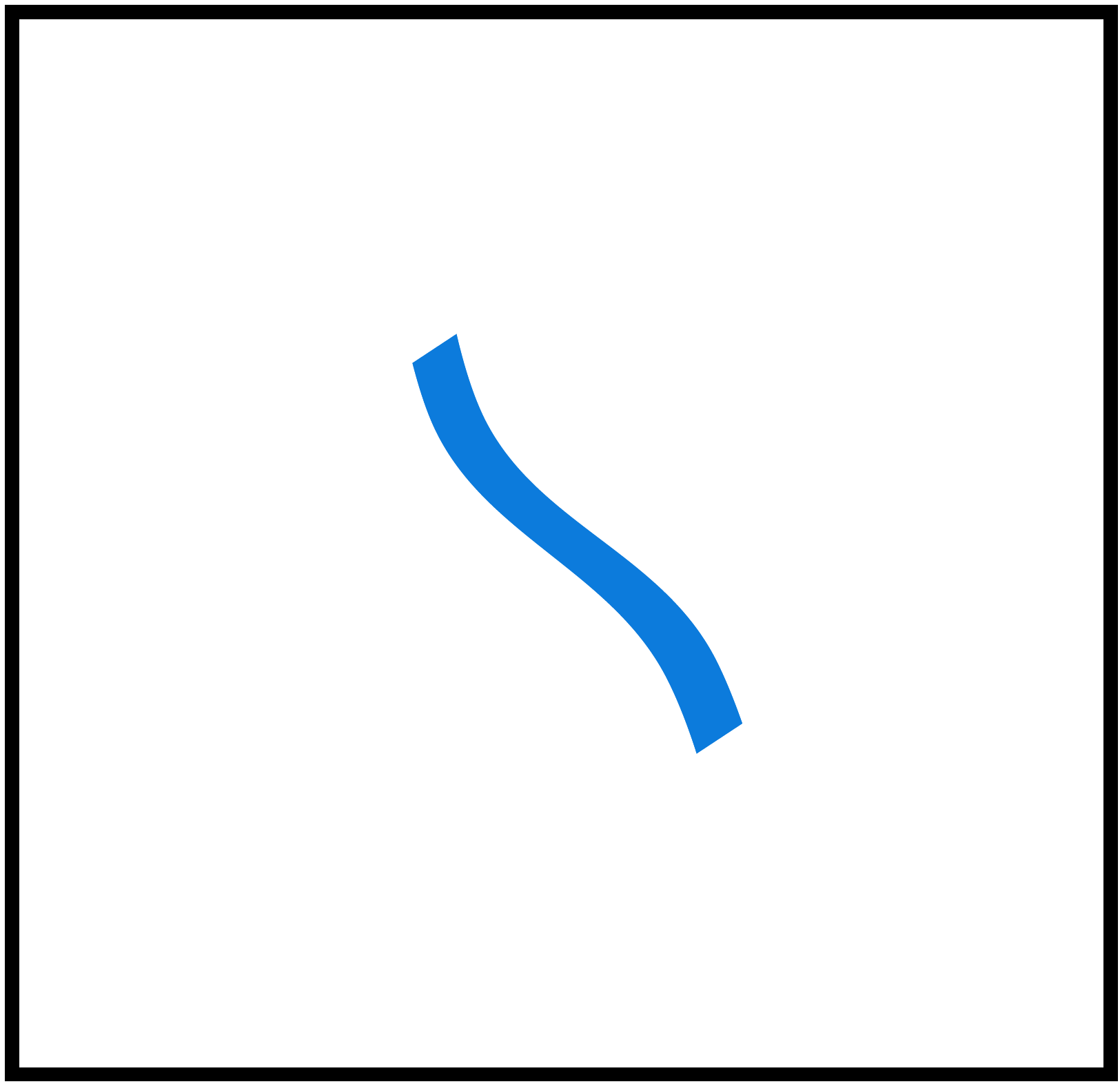}} \quad
    \caption{
    Visualization of the spaND algorithm. 
    Interfaces are shown in blue and interiors in yellow.
    At each level, interiors only interact with (i.e., have edges connecting to) the neighboring interfaces.
    Interiors are first eliminated using the block Gaussian elimination.
    The remaining interfaces are then sparsified, which effectively eliminates (without fill-in) a number of
        variables from the system.
    The algorithm then proceeds to the next level and completes when the last level interface is eliminated.
    }
    \label{fig:spand-illustrate}
\end{figure}

We now describe the spaND algorithm, the details of which can be found in \cite{cambier2020algebraic}.
The algorithm is illustrated in \Cref{fig:spand-illustrate}.
It uses a multilevel ordering based on nested dissection. 
This ordering (called the \emph{spaND partitioning}) is fully algebraic, that is, only the entries of $A$ are needed to define it (more specifically, the partitioning is defined on the adjacency graph of $A$).
At each level, the unknowns are split into two groups of subsets: interiors and interfaces.
The interfaces are small subsets that separate the interiors from each other.
The precise definitions of interiors and interfaces can be found in \cite{cambier2020algebraic}.

Once the multilevel partitioning has been defined, one can perform the actual factorization, see \Cref{fig:spand-illustrate}.
At each level, the interiors are first eliminated using the block Gaussian elimination. 
The existence of interfaces limits the fill-in (Schur complement updates) arising during elimination. 
The fill-in matrix blocks (interactions between interfaces) 
are scaled as in \cref{eq:prescale}, and afterward sparsified using the approximation scheme of \cref{sec:main}, which  reduces the size of the interfaces.
At this point the algorithm proceeds to the next level.
The algorithm completes when the variables in the last level are eliminated using the Cholesky decomposition.
The repeated reduction in the sizes of the interfaces 
allows obtaining an efficient sparse algorithm.
 


\begin{figure}
    \centering
    \begin{tikzpicture}
    \node (image) at (0,0) {\includegraphics[width=0.3\textwidth]{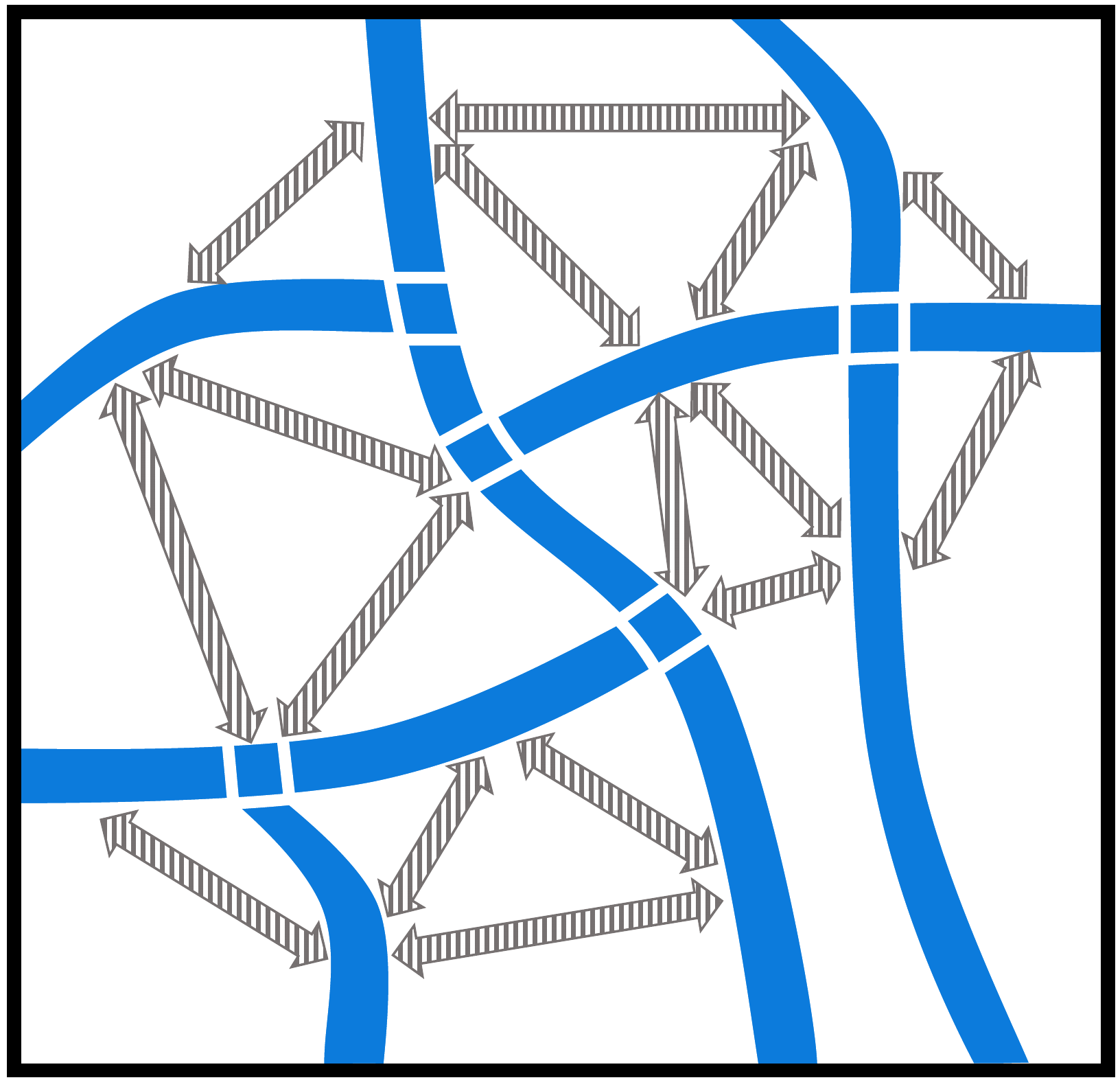}};
    \node at (-2.7,1.7) {$p$};
    \draw[thick,-latex] (-2.5,1.7) to (-0.6,1.7);
    \node at (-2.7,1.2) {$A_{pw}$};
    \draw[thick,-latex] (-2.3,1.1) to (-1.15,1.0);
    \draw[thick,-latex] (-2.3,1.2) to (-0.1,1.1);
    \draw[thick,-latex] (-2.3,1.3) to (0.0,1.48);
    \end{tikzpicture}
    \begin{tikzpicture}
    \node (image) at (0,0) {\includegraphics[width=0.3\textwidth]{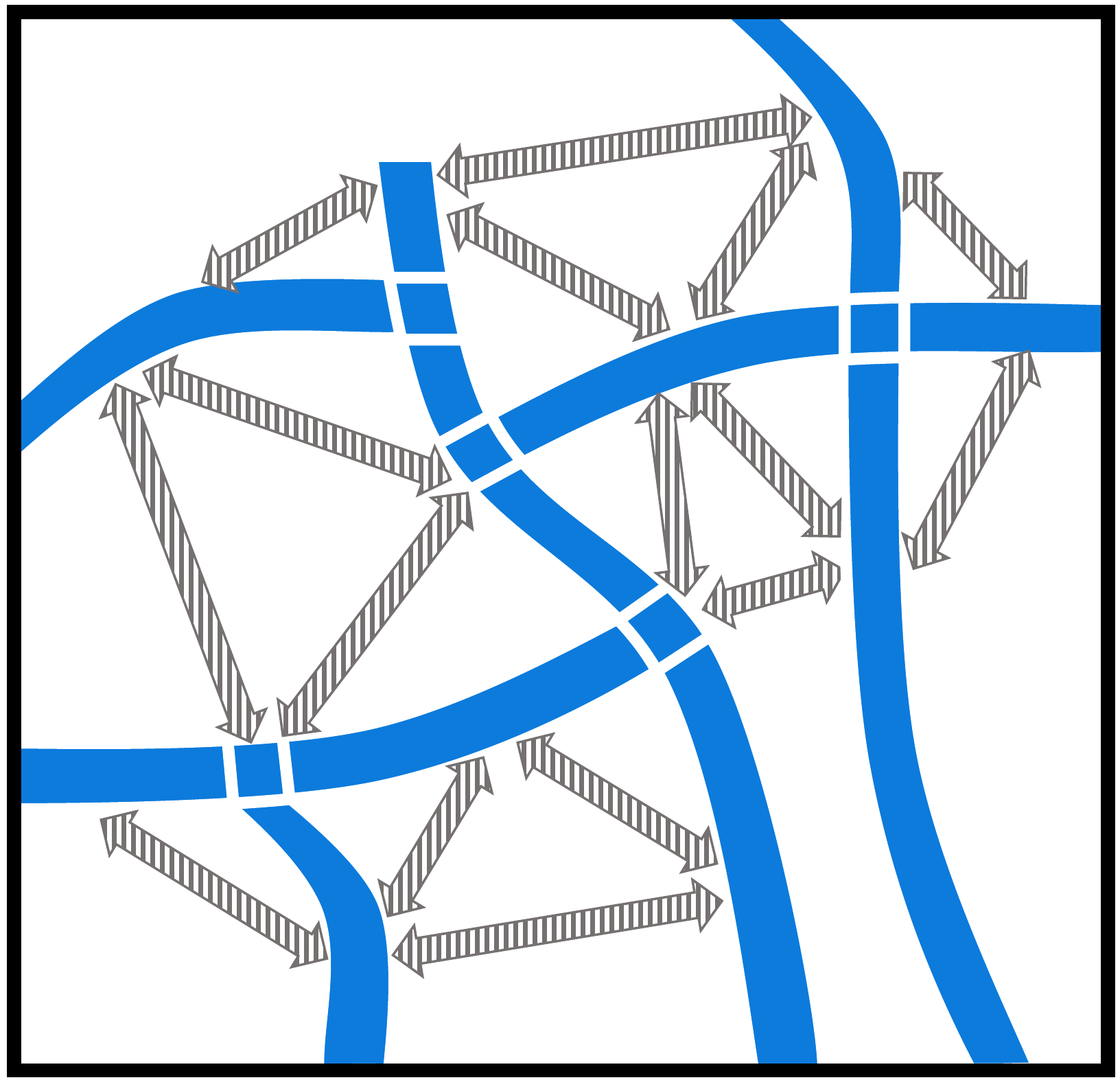}};
    \node at (2.7,1.1) {$c$};
    \draw[thick,-latex] (2.5,1.1) to (-0.5,1.1);
    \end{tikzpicture}
    \caption{
    	Illustration of the sparsification of a single interface $p$.
        The interactions between interfaces are conceptually shown in gray. The block
        $A_{pw}$ denotes the interactions of $p$ with the neighboring interfaces.
        Sparsifying $p$ splits it into $c$ (coarse variables) and $f$ (fine variables).
        Since $f$ is disconnected from all other interfaces, it is effectively eliminated 
        ($f$ is represented by the piece of $p$ dropped in the picture above).
        This reduces the size of the interface without introducing any fill-in.
        }
        \label{fig:spand-sparsify}
\end{figure}
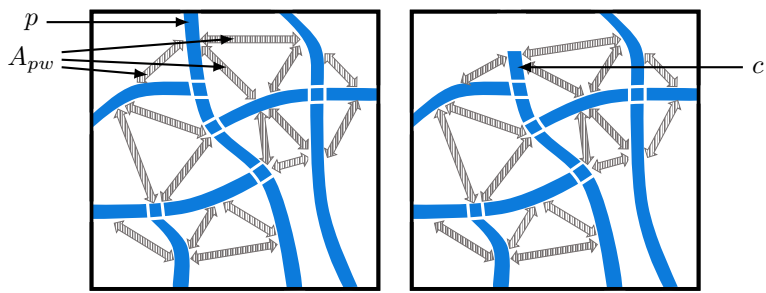

The sparsification process is depicted in \cref{fig:spand-sparsify}.
After eliminating interiors, most of connections between interfaces are through fill-in blocks. 
Consider the interface $p$ in \cref{fig:spand-sparsify}.
The connections to its neighbors (collectively denoted by $w$) typically have the low-rank property.
Using notation from \cref{sec:main} this means that the block $A_{12} := A_{pw}$ has a quickly decaying spectrum.
We can therefore sparsify $p$ as described in \cref{sec:main}. 
The step described in \cref{eq:Q-basis} is a change of basis that splits the variables from $p$ into $c$ (the coarse variables) and $f$ (the fine variables). 
After the approximation step ---\cref{eq:trailing-1} (first-order), \cref{eq:trailing-2} (full second-order), or \cref{eq:trailing-3} (superfine second-order)--- has been applied, the variables from $f$ are disconnected from all other unknowns and effectively eliminated.
This process does not introduce any fill-in.
\Cref{algo:spand} describes the spaND algorithm using a recursive formulation.


\subsection{Obtaining the $E$ or $E_2$ matrix}
To compute the $Q$ matrix from \cref{eq:Q-basis}, spaND uses the column-pivoted rank-revealing QR (RRQR) which gives $A_{12}\Pi = QR$ where $\Pi$ is a permutation matrix. When using the first-order or the \fullsec, the decomposition can be stopped after $k$ steps, where $k$ is the smallest number such that $|R(k,k)| < \varepsilon |R(1,1)|$ (such relative criterion is typically chosen \cite{cambier2020algebraic,klockiewicz2019sparse,ho2016hierarchical,pouransari2017fast,chen2018distributed,feliu2018recursively}). At this point, we have
\begin{equation}
	A_{12} \Pi = 
		\begin{pmatrix}
		Q_c & Q_{f}
	\end{pmatrix}
	R =
	\begin{pmatrix}
		Q_c & Q_{f}
	\end{pmatrix}
		\begin{pmatrix}
		R_{cc}     & R_{cf}	 	\\
	    0       		& R_{f f}	\\
	\end{pmatrix} 
\end{equation}
where $Q_c$ has $k$ columns, and 
\[ E =  Q_f^\top A_{12} = \begin{pmatrix} 0 & R_{ff} \end{pmatrix} \Pi^\top\] 
Since both $R_{ff}$ and $\Pi$ are computed in the RRQR routine (see \Cref{sec:experiments}), obtaining $E$ involves a negligible additional cost.
We also have $Q_c^\top A_{pw} = \begin{pmatrix} R_{cc} & R_{cf} \end{pmatrix} \Pi^\top$.

In the \superfine, the decomposition is stopped later, with $k$ being the smallest number such that $|R(k,k)| < \varepsilon^2 |R(1,1)|.$ At this point
\begin{equation}\label{eq:rrqr}
	A_{12} \Pi = 
		\begin{pmatrix}
		Q_c & Q_{f_2} & Q_{f_1}
	\end{pmatrix}
	R =
	\begin{pmatrix}
		Q_c & Q_{f_2} & Q_{f_1}
	\end{pmatrix}
		\begin{pmatrix}
		R_{cc}     & R_{cf_2}	 	& R_{cf_1}\\
	    0       		& R_{f_2 f_2}	& R_{f_2 f_1}\\
	    0				&	0					& R_{f_1 f_1}
	\end{pmatrix} 
\end{equation} We only need
\[
E_2 	 = Q_{f_2}^\top A_{12} = 
\begin{pmatrix}
    0 	&  R_{f_2 f_2} & R_{f_2 f_1} \\
\end{pmatrix} \Pi^\top
\] 
Should $Q$ be computed separately, e.g., using randomized SVD, the matrix $E$ (or $E_2$) is obtained by performing the multiplication $Q_f^\top A_{12}$ (or $Q_{f_2}^\top A_{12}$). 
The (typically larger) multiplication $Q_c^\top A_{12}$ has to be performed regardless, and in our experience, is a  small portion of the computations.
Thus the second-order scheme involves only a minor additional computational cost also in this case.

\begin{algorithm}
    \begin{algorithmic}
    \REQUIRE $A_{pw}, 0 \leq \varepsilon \leq 1$
                    \[
            A_{pw} \Pi = \begin{cases}
            		\begin{pmatrix}
		Q_c & Q_{f}
	\end{pmatrix}
		\begin{pmatrix}
		R_{cc}     & R_{cf}	 	\\
	    0       		& R_{f f}	\\
	\end{pmatrix}  & \parbox[t]{12em}{for the \textbf{\fo} \\ or the  \textbf{\fullsec}} \\
		\begin{pmatrix}
		Q_c & Q_{f_2} & Q_{f_1}
	\end{pmatrix}
		\begin{pmatrix}
		R_{cc}     & R_{cf_2}	 	& R_{cf_1}\\
	    0       		& R_{f_2 f_2}	& R_{f_2 f_1}\\
	    0				&	0					& R_{f_1 f_1}
	\end{pmatrix} & \parbox[t]{12em}{for the \textbf{\superfine}}
            		\end{cases}
            \]
     where $\| R_{ff} \| = \mathcal{O}(\varepsilon)$, $\| R_{f_1 f_1} \| = \mathcal{O}(\varepsilon^2)$. 
            \[
            \tE = \begin{cases}
            		0 						   & \mbox{for the \textbf{\fo}} \\
            		E = Q_f^\top A_{pw} = \begin{pmatrix} 0 & R_{ff} \end{pmatrix} \Pi^\top  & \parbox[t]{12em}{for the \textbf{\fullsec}} \\
            		
            		\begin{pmatrix} 0 \\ E_2 = Q_{f_2}^\top A_{pw} \end{pmatrix} =
            		\begin{pmatrix} 0 \\
    0 	&  R_{f_2 f_2} & R_{f_2 f_1} \\
\end{pmatrix} \Pi^\top & \parbox[t]{12em}{for the \textbf{\superfine}}
            		\end{cases}
            \]
    \STATE $Q_c^\top A_{pw} = \begin{pmatrix} R_{cc} & R_{cf} \end{pmatrix} \Pi^\top$
    \RETURN $Q_c, Q_f, \tE, Q_c^\top A_{pw}$
    \end{algorithmic}
    \caption{ Run column-pivoted QR; compute $E$ (or $E_2$), $Q_c$, $Q_f$ and $Q_c^\top A_{pw}$ }
    \label{algo:sparsify}
\end{algorithm}


\subsection{Accuracy and relation between the factorization and solve}
Notice that $\varepsilon > 0$ controls the accuracy of the low-rank approximations and of the entire algorithm.
Smaller value results in a more accurate factorization, which will take more time to compute. The resulting preconditioner will be more expensive to apply, but it will approximate $A^{-1}$ more accurately, resulting in fewer iterations in the solve phase (e.g., using Conjugate Gradient).
We have $E \to 0$ as $\varepsilon \to 0$, in which case spaND becomes an exact block Cholesky factorization with a nested dissection ordering, and gives the exact solution.
The optimal total runtime will likely be obtained long before, however, when the runtimes of the factorization and solve phases are more balanced.

\begin{algorithm}
    \begin{algorithmic}
    \REQUIRE $A$, $0 \leq\varepsilon \leq 1$
        \STATE \textbf{spaND}($A^l$, $l$)
        \FORALL[$\texttt{Eliminate interiors}$]{$s$ interior}
        \STATE Eliminate $s$ in one step of block Cholesky algorithm (with $A_{ss} = L_{s} L_{s}^\top$)
            \[ \begin{matrix}
            	G_s \begin{pmatrix} A_{ss} & A_{sw} \\ A_{ws} & A_{ww} \end{pmatrix} G_s^\top = \begin{pmatrix} I & \\ & A_{ww} - A_{ws} A_{ss}^{-1} A_{sw} \end{pmatrix}, & G_s = \begin{pmatrix} L_{s}^{-1} & \\ - A_{ws} L_{s}^{-\top} & I \end{pmatrix} 
            	\end{matrix}
            \]
        \ENDFOR
        \FORALL[\texttt{Scale interfaces}]{$p$ interface}
        \STATE Scale $p$ using the Cholesky algorithm (with $A_{pp} = Z_{p} Z_{p}^\top$)
            \[ \begin{matrix}
                	D_p \begin{pmatrix} A_{pp} & A_{pw} \\ A_{wp} & A_{ww} \end{pmatrix} D_p^\top = \begin{pmatrix} I & Z_{p}^{-1} A_{pw} \\ A_{wp} Z_{p}^{-\top} & A_{ww}\end{pmatrix}, & D_p = \begin{pmatrix} Z_{p}^{-1} & \\ & I \end{pmatrix}        	
            	\end{matrix}
			 \]
        \ENDFOR
        \FORALL[\texttt{Sparsify interfaces}]{$p$ interface}
        \STATE Run rank-revealing QR (\Cref{algo:sparsify}) to sparsify $p$ 
            \[ E_p Q_p^\top \begin{pmatrix} I & A_{pw} \\ A_{wp} & A_{ww} \end{pmatrix} Q_p E_p^\top \approx \begin{pmatrix} I & & \\ & I & Q_{c}^\top A_{pw}\\ & A_{wp} Q_{c} & A_{ww} \end{pmatrix}  \]
            \[ \begin{matrix}
            	Q_p = \begin{pmatrix} \begin{pmatrix} Q_f & Q_c \end{pmatrix} & \\ & I \end{pmatrix}, 
            	& E_p = \begin{pmatrix} I & & \\ & I & \\ - \tE^\top & & I \end{pmatrix} 
            	\end{matrix} \]         	
        \ENDFOR
        \STATE 
        $	\begin{matrix}
				A^{l-1} =  B_l A^l B_l^\top, & B_l = \prod_p E_p^l Q_p^{l,\top} \prod_p D_p^l \prod_s G_s^l
			\end{matrix} $
        \IF{$l > 0$}
        \STATE Recurse \textbf{spaND}($A^{l-1}_{P,P}$, $l-1$) where $P$ denotes the rows/columns corresponding to the non-yet-eliminated interfaces
        \ENDIF
    \end{algorithmic}
    \caption{
    The spaND algorithm in a recursive form.
    We assume that the special spaND partitioning \cite{cambier2020algebraic} on $\ell \geq 1$ levels has been computed. 
    The algorithm starts at level $l = \ell$, with $A^\ell = A$, and completes at $l = 0$. 
    The result is $M = L^{-\top}L^{-1}$ such that $L ^{-1}A L^{-\top} \approx I$ with $L^{-1} = \prod_{l=0}^{\ell} \left( \prod_p (E_p^l Q_p^{\top,l}) \prod_p D_p^l \prod_s G_s^l \right)$}
    \label{algo:spand}
\end{algorithm}

%% file: experiments.tex
\section{Experimental results}
\label{sec:experiments}
We compare the preconditioners obtained when using the first- and second-order approximation schemes in the spaND algorithm. In all tested cases, the number of levels $\ell$ in the spaND partitioning is the closest integer to $\log_2(n/25),$ where $n$ is the number of rows in the tested matrix.
We also skip the scaling and sparsification of interfaces in the first four levels of the algorithm, when the interfaces are still small.
Thus, the only varying parameter is the accuracy $0 \leq \varepsilon \leq 1$.

The spaND implementation is sequential and was written in C++, using
BLAS and LAPACK \cite{anderson1999lapack} routines provided by Intel(R) MKL.
In particular, the (early-stopping) rank-revealing QR factorization is implemented using the $\texttt{dlaqps}$ routine from LAPACK.
All experiments were run on CPUs with Intel(R) Xeon(R) E5-2640v4 (2.4GHz) processor with 128 GB RAM, always using a single thread.

The approximate inverse operator $M = L^{-\top}L^{-1}$ returned by spaND, is used in the preconditioned Conjugate Gradient (PCG) with a zero initial guess and convergence declared when the relative 2-norm of the residual falls below $10^{-10}$.

We further use the following notation:
\begin{itemize}
	\item $n$ -- is the number of unknowns (rows) of the given matrix $A$
	\item $nnz(A)$ -- is the number of nonzero entries of $A$
	\item $n_{cg}$ -- denotes the number of PCG iterations needed to converge
	\item $\mu$ -- denotes the memory requirements needed to store the preconditioner, defined as $\mu := nnz(L)/nnz(A)$
	\item $T_f$ -- is the time needed to perform the spaND hierarchical factorization
	\item $T_s$ -- is the time needed by PCG to converge
	\item $T_{t}$ -- is the total time needed to solve the system, i.e., $T_t = T_f + T_s$
	\item $\Delta T_s$ -- is the difference (in \%) in the total runtime, when compared to the \fo
 \end{itemize}
All times are reported in seconds. 

\subsection{Low- and high-contrast Laplacians}
\begin{figure}[htbp]
\centering
\includegraphics[width=0.2\textwidth,height=0.2\textwidth]{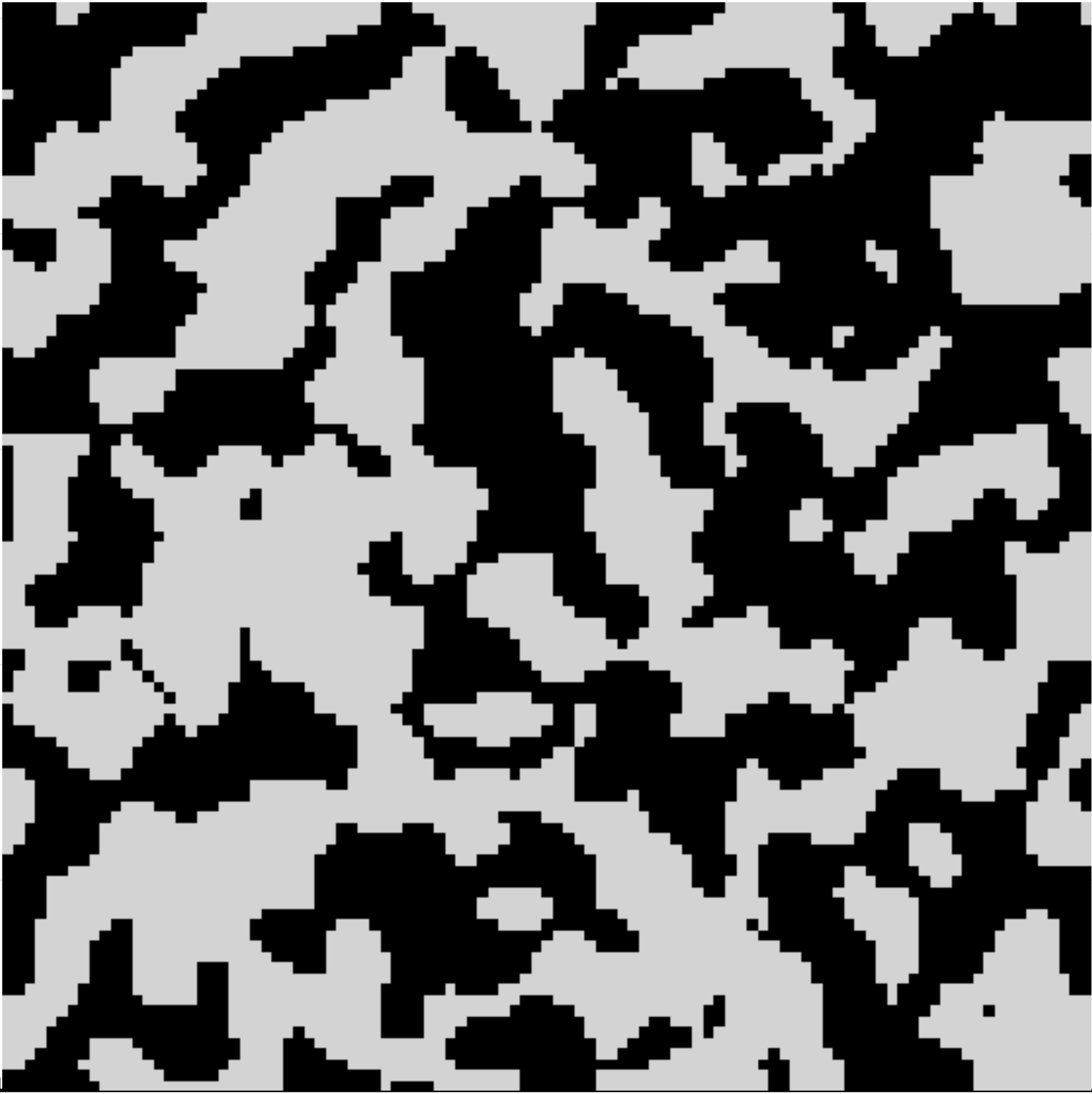}
\includegraphics[width=0.2\textwidth,height=0.2\textwidth]{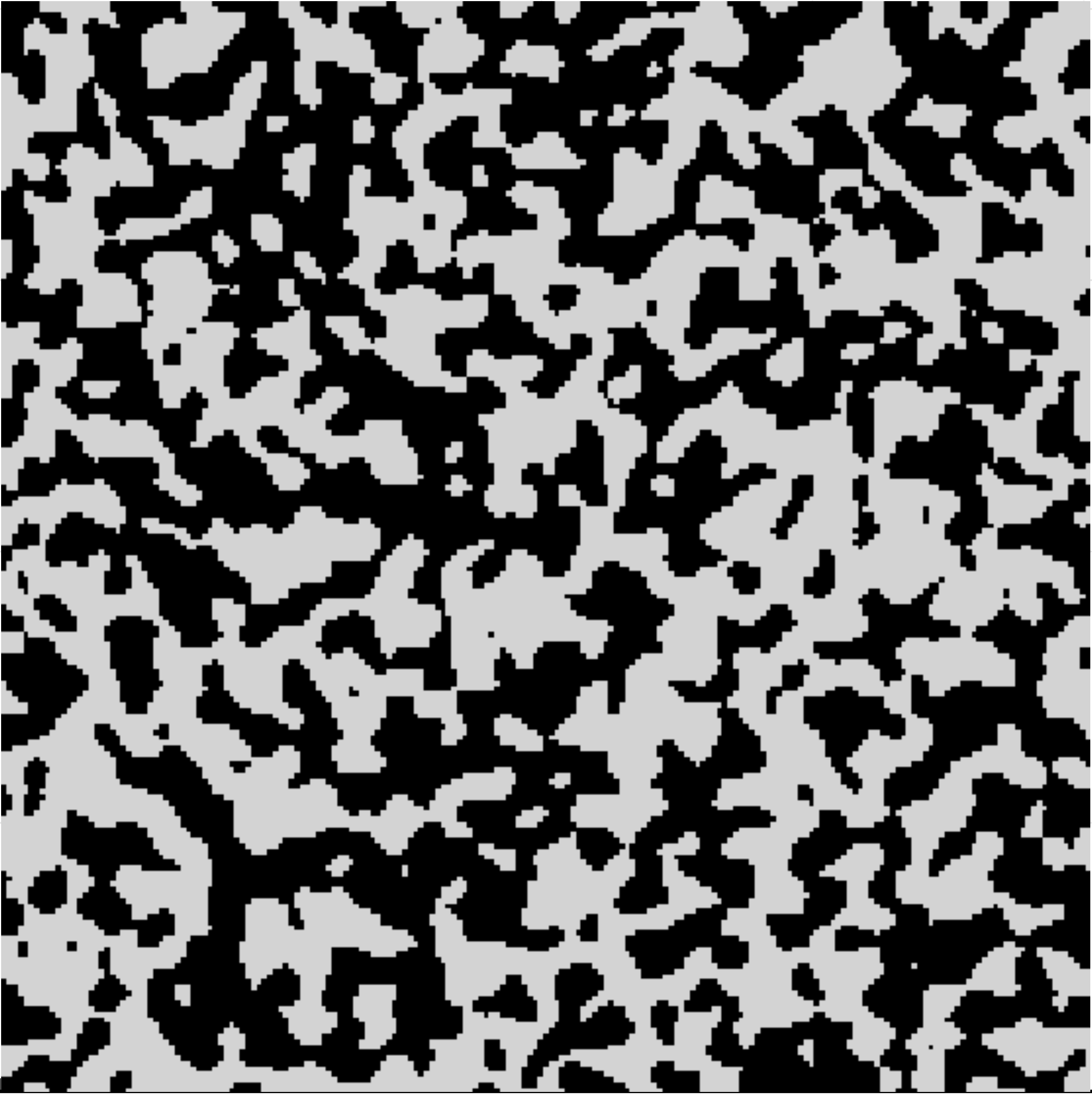}
\caption{Examples of random coefficient fields used to define high-contrast Laplacians for $d = 50$ (left) and $d = 100$ (right). }
\label{fig:high-contrast}	
\end{figure}

We consider the elliptic equation
\begin{equation}
\begin{matrix} \nabla (a(x) \cdot \nabla u(x)) = f	 & \forall{x \in \Omega \in \left[ 0,1\right]^2}, & u|_{\partial \Omega} = 0
\end{matrix}	
\end{equation}
discretized using the 5-point stencil method, on a square $d \times d$ grid.
For $\rho \geq 1.0$, we define the $a(x)$ field as in \cite{cambier2020algebraic,feliu2018recursively}. Namely, for $i,j = 1,2, \ldots, d,$ we pick a random $\hat{a}_{ij} \in (0,1)$ from the uniform distribution, and convolve the resulting $d \times d$ field $\hat{a}$ with an isotropic Gaussian of standard deviation $\sigma= 2$ to smooth the field out. We then quantize
\begin{equation}
a_{ij} =
	\begin{cases}
		\rho & \mbox{if } \hat{a}_{ij} \geq 0.5 \\
		\rho^{-1} & \mbox{if } \hat{a}_{ij} < 0.5.
	\end{cases}
\end{equation}
Thus $\rho$ is a contrast parameter of the field. At $\rho = 1.0$ we obtain the constant coefficient Poisson equation since then $a(x) \equiv 1.0$. At $\rho = 100.0$ the contrast between coefficients is $\rho^2 = 10^4$. Examples of the coefficient fields are shown in \cref{fig:high-contrast}. The condition number $\kappa(A)$ of the resulting matrix $A$ scales approximately as $ \rho^2 n.$

\subsubsection{Forward errors on eigenvectors}
For the constant-coefficient case, i.e., when $\rho=1.0$, the unit-length eigenvectors of $A$ are known exactly. We therefore compute the forward errors $ \| (I - M A) v_{\lambda} \|_2$ on selected unit-length eigenvectors $v_\lambda$. More precisely, for a given $p \in \{1,2,\ldots,d\}$ we consider the eigenvector corresponding to the function of the grid given by
\[
	x(i,j) = \sin \left(\frac{pi\pi}{d+1}\right) \sin\left(\frac{pj\pi}{d+1} \right)
\]
which we normalize to obtain a unit-norm eigenvector $v_{\lambda_p}$. The corresponding eigenvalue is $\lambda_p = 8 \sin^2\left(\frac{p \pi}{2(d+1)}\right).$
In \Cref{fig:poisson-fwd-error} we plot the forward error as a function of the corresponding eigenvalue $\lambda$ for $d = 1000$, i.e., for the $1000 \times 1000$ grid.
Because computing errors on all $10^6$ eigenvectors would be infeasible and difficult to plot, we consider $p \in \{ \lfloor (5/4)^k \rfloor : k = 0,1,2, \ldots \}, \ p \leq d$.
The corresponding eigenvalues fall in the whole range of magnitudes. 
For a given value of $\varepsilon$ parameter, the (full) second-order scheme is more accurate on all of the spectrum compared to the first-order scheme.
The difference is particularly pronounced on the middle-to-high frequency eigenmodes.
The accuracy on the lowest-frequency eigenmodes depends largely on the accuracy parameter $\varepsilon$.


\pgfplotstableread[col sep=comma]{logs/poisson-0.1-FO-noprescale.csv}\datapf
\pgfplotstableread[col sep=comma]{logs/poisson-0.1-SO-noprescale.csv}\dataps
\pgfplotstableread[col sep=comma]{logs/poisson-0.01-FO-noprescale.csv}\datapfff
\pgfplotstableread[col sep=comma]{logs/poisson-0.01-SO-noprescale.csv}\datapsss
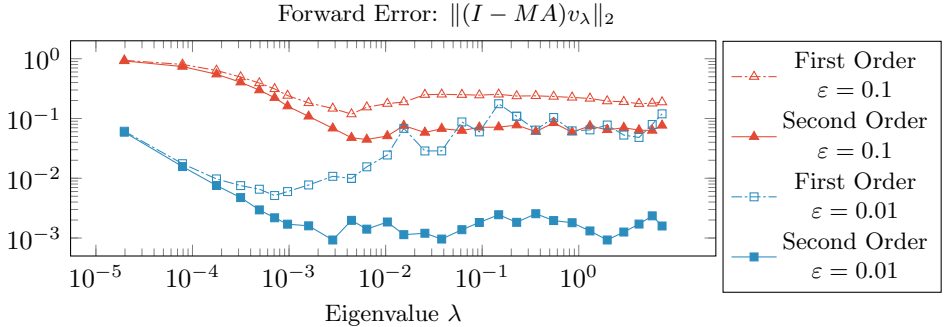
\begin{figure}[htbp]
\centering
	\begin{tikzpicture}
	\begin{axis}[
		width=4in,
		height=1.75in,
		ymin=5*10^(-4),
		ymode=log,
		xmode=log,
		xtick={10^(-5),10^(-4),10^(-3),10^(-2),10^(-1),10^0},
		ytick={10^(-3),10^(-2),10^(-1),10^0},
		legend pos=north east,
		legend entries={\small First Order \\ $\varepsilon = 0.1$,
								 \small Second Order \\ $\varepsilon = 0.1$,
								 \small First Order \\ $\varepsilon = 0.01$,
								 \small Second Order \\ $\varepsilon = 0.01$},
		legend style={font=\small, align=center, at={(1.35,1.0)}},
		xlabel={\small Eigenvalue $\lambda$},
	]
	\addplot[eps01fo]	table[x=axnorm, y=eigfwerr]{\datapf};
	\addplot[eps01so]	table[x=axnorm, y=eigfwerr]{\dataps};
	\addplot[eps001fo]	table[x=axnorm, y=eigfwerr]{\datapfff};
	\addplot[eps001so]	table[x=axnorm, y=eigfwerr]{\datapsss};
	\end{axis}
	\node(title) at (5.0,3.2) {\small Forward Error: {$\|(I - MA) v_\lambda \|_2$}};
	\end{tikzpicture}
	\caption{Forward errors on the unit-length eigenvectors of the 2D constant-coefficient Poisson equation. On the middle-to-high-frequency eigenmodes, the (full) second-order scheme with $\varepsilon = 0.1$ is equally accurate as the first-order scheme with $\varepsilon = 0.01$.}
	\label{fig:poisson-fwd-error}
\end{figure}

\subsubsection{Halved PCG iteration counts and improved total timings}
We perform a scaling study on square grids of increasing sizes for $\rho = 1.0$ (constant-coefficient field) and $\rho = 100.0$ (high-contrast field).
The problems are refined by doubling the size of the grid in each dimension.
In each case we solve the system $Ax = b$ where $b$ is a vector of ones.
The results for $\varepsilon = 0.01$ and $\varepsilon = 0.001$ are shown in \cref{tab:highcontrast}.
For both values of $\rho$, and all tested grid sizes, we consistently observe that the number of iterations needed for convergence is halved when using the \fullsec{} as compared to the \fo, with approximately the same factorization time, resulting in improved total timings.

\bgroup
\setlength{\tabcolsep}{0.35em}
\begin{table}[htbp]
\renewcommand{\arraystretch}{1.1}
\small
\centering {
	\begin{tabular}{ccc | rrrrr | rrrrrr}
	  \multicolumn{5}{c}{$\bf \rho=1.0$ \textbf{(no contrast)}} &\multicolumn{3}{c}{First Order} & \multicolumn{5}{c}{Second Order} \\
	 $d$ &  $n = d^2$ & $\varepsilon$ & $\mu$ & $n_{cg}$ & $T_f[s]$ & $T_s[s]$ & $T_t[s]$ & $\mu$ & $n_{cg}$ & $T_f[s]$ & $T_s[s]$ & $T_{t}[s]$ \\
	 \hline
	400&0.16M&0.01&7.8&\textbf{9}&0.5&0.6&1.1&8.6&\textbf{5}&0.4&0.3&0.6\\
	800&0.64M&0.01&7.7&\textbf{11}&2.1&2.9&5.1&8.5&\textbf{6}&1.7&1.5&3.2\\
	1600&2.56M&0.01&7.7&\textbf{16}&8.8&17.4&26.2&8.5&\textbf{8}&10.1&10.6&20.7\\
	3200&10.2M&0.01&7.7&\textbf{22}&37.1&103.9&141.0&8.5&\textbf{11}&35.9&54.8&90.7\\
	6400&41.0M&0.01&7.6&\textbf{34}&145.0&666.9&811.9&8.4&\textbf{17}&150.0&346.1&496.1\\
	\hline
	400&0.16M&0.001&8.1&\textbf{5}&0.5&0.3&0.8&8.9&\textbf{3}&0.4&0.2&0.6\\
	800&0.64M&0.001&8.0&\textbf{6}&1.8&1.4&3.1&8.8&\textbf{3}&2.0&0.9&2.9\\
	1600&2.56M&0.001&8.0&\textbf{7}&7.5&6.3&13.8&8.9&\textbf{4}&9.3&4.9&14.2\\
	3200&10.2M&0.001&8.0&\textbf{8}&35.6&35.8&71.4&8.8&\textbf{4}&35.8&19.9&55.8\\
	6400&41.0M&0.001&7.9&\textbf{10}&167.6&215.7&383.3&8.7&\textbf{5}&152.4&115.8&268.2\\
	 \hline \hline
	  \multicolumn{5}{c}{$\bf \rho=100.0$ \textbf{(high contrast)}} &\multicolumn{3}{c}{First Order} & \multicolumn{5}{c}{Second Order} \\
	 $d$ &  $n = d^2$ & $\varepsilon$ & $\mu$ & $n_{cg}$ & $T_f[s]$ & $T_s[s]$ & $T_t[s]$ & $\mu$ & $n_{cg}$ & $T_f[s]$ & $T_s[s]$ & $T_{t}[s]$ \\
	 	 \hline
	400&0.16M&0.01&7.6&\textbf{15}&0.5&0.9&1.4&8.3&\textbf{7}&0.5&0.4&0.9\\
	800&0.64M&0.01&7.5&\textbf{22}&2.1&5.7&7.8&8.3&\textbf{11}&1.8&2.9&4.7\\
	1600&2.56M&0.01&7.6&\textbf{28}&8.7&28.2&36.9&8.3&\textbf{13}&7.6&13.7&21.3\\
	3200&10.2M&0.01&7.5&\textbf{46}&33.6&197.0&230.6&8.3&\textbf{22}&29.6&90.4&120.0\\
	6400&41.0M&0.01&7.5&\textbf{82}&140.6&1599.0&1739.6&8.2&\textbf{38}&142.6&748.5&891.1\\ 
	\hline	 	 
	400&0.16M&0.001&7.8&\textbf{8}&0.4&0.4&0.8&8.5&\textbf{4}&0.5&0.3&0.7\\
	800&0.64M&0.001&7.7&\textbf{9}&1.8&2.2&4.0&8.5&\textbf{5}&1.8&1.2&3.0\\
	1600&2.56M&0.001&7.8&\textbf{10}&7.5&9.1&16.6&8.5&\textbf{5}&7.4&5.0&12.5\\
	3200&10.2M&0.001&7.7&\textbf{12}&38.4&53.8&92.2&8.5&\textbf{6}&35.3&29.4&64.7\\
	6400&41.0M&0.001&7.7&\textbf{16}&144.9&316.4&461.3&8.5&\textbf{8}&154.7&182.9&337.6\\
\end{tabular}
	}
\caption{Scaling study on 2D Laplacians. Second-order scheme consistently halves the number of PCG iterations for approximately the same factorization time.}
\label{tab:highcontrast}
\end{table}
\egroup

\subsection{SuiteSparse matrices}
 To test the efficiencies of the new approximation schemes in practice, we run the spaND algorithm on all SPD matrices from the University of Florida sparse matrix collection \cite{davis2011university} (SuiteSparse), with at least 500,000 rows.
We run spaND with the \fo, the \fullsec{}, and the \superfine. 
We test four values of the accuracy parameter $\varepsilon=0.2,0.1,0.05,0.01$.
At $\varepsilon = 0.01,$ PCG converges in a small number of iterations on all tested matrices.
As mentioned above, $\varepsilon$ is the only varying parameter. 

In each case, we solve the system $A' x = b'$ where $A' = D^{-\frac{1}{2}} A D^{-\frac{1}{2}}$, $b' = D^{-\frac{1}{2}}b$ with $D$ being the diagonal of $A$, and  $b$ a vector of ones. 
Such diagonal prescaling is recommended 
\cite{vanek1999two,franceschini2019block}, when solving structural problems which are among the most challenging test cases.

\subsubsection{Halved PCG counts and improved solve times}
Similar as in the Laplace scaling study, for a given accuracy parameter $\varepsilon,$ the number of PCG iterations needed for convergence is almost exactly halved when using the \fullsec, across all tested cases.
This is also true for the \superfine, on almost all of the problems.
Example aggregate results are shown in \Cref{fig:cg-suite-relative}.

For a given $\varepsilon$ parameter the preconditioners using the new second-order schemes are more expensive to apply.
However, the increase in memory requirement (proportional to the cost of applying the preconditioner), is moderate, never exceeding 100\%, as shown in \Cref{fig:density-suite-relative}.
As a result, the time needed for convergence of PCG is still significantly reduced in all cases.
This is shown in \Cref{fig:solve-suite-relative}.

\subsubsection{Squeezed convergence plots}
In \Cref{fig:cg-convergence} we show example plots of the 2-norm error decay in the PCG algorithm. 
We observe that, for a given problem and accuracy parameter, the behavior of the residual $\| A' x_k - b' \| / \| b' \|$ as a function of $k$, when using the \fullsec{}, is approximately the same as the behavior of $\| A' x_{2k} - b' \| / \| b' \| $ when using the \fo{}.
As a result, the plot when using the \fullsec{} is squeezed compared to the plot when using the \fo{}, but retains its general shape.
This behavior is approximately the same as in the conclusion of \cref{thm:residuals}.
We note here that we observed the same behavior of convergence plots when the right-hand sides were chosen at random.

\pgfplotstableread[col sep=comma]{logs/suitesparse-3-0.1-prescale-SO-best-relative.csv}\datax
\pgfplotstableread[col sep=comma]{logs/suitesparse-3-sf-0.1-prescale-SO-best-relative.csv}\dataxx
\pgfplotstableread[col sep=comma]{logs/suitesparse-3-0.05-prescale-SO-best-relative.csv}\datay
\pgfplotstableread[col sep=comma]{logs/suitesparse-3-sf-0.05-prescale-SO-best-relative.csv}\datayy
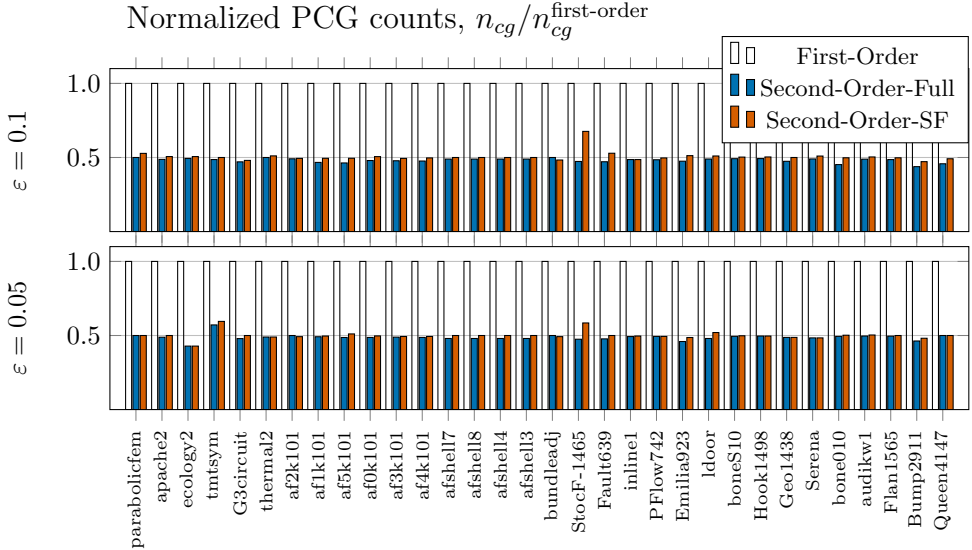
\begin{figure}[htbp]
	\begin{tikzpicture}
	\begin{groupplot}[
		group style={
        group name=suite-sparse-1,
        group size=1 by 2,
        xlabels at=edge bottom,
        xticklabels at=edge bottom,
        vertical sep=0cm,
        ylabels at=edge left,
        yticklabels at=edge left,
        vertical sep=0.2cm
    },
    xmin=-1,
    xmax=32,
    ymin=0,
    ymax=1.1,
    ybar=0.5pt,
    ytick={0.5,1.0},
    yticklabels={0.5,1.0},
    yticklabel style={rotate=-90},
    width=0.95\textwidth,
    xtick=data,
    xticklabels from table={\datay}{matrix},
    xticklabel style={tick label style={rotate=90}, font=\scriptsize},
    xtick={0,1,2,3,4,5,6,7,8,9,10,11,12,13,14,15,16,17,18,19,20,21,22,23,24,25,26,27,28,29,30,31},
    ]
    \nextgroupplot[
    	bar width=0.08cm,
    	height=3.75cm,
    	width=13.0cm,
    	legend entries={First-Order, Second-Order-Full, Second-Order-SF},
    	legend style={at={(1.0,1.2)},
    							font=\small},
    	ylabel style={align=left},
    	ylabel={$\varepsilon=0.1$},
    	ymajorgrids,
    ]
	\addplot[color=black,fill=none]	table[x expr=\coordindex, y expr=1.0]{\datax};
	\addplot[color=black,fill=colorz]	table[x expr=\coordindex, y=cg_rel]{\datax};
	\addplot[color=black,fill=colorx]	table[x expr=\coordindex, y=cg_rel]{\dataxx};
	    \nextgroupplot[
    	bar width=0.08cm,
    	height=3.75cm,
    	width=13.0cm,
    	ylabel style={align=left},
    	ylabel={$\varepsilon=0.05$},
    	ymajorgrids,
    ]
	\addplot[color=black,fill=none]	table[x expr=\coordindex, y expr=1.0]{\datay};
	\addplot[color=black,fill=colorz]	table[x expr=\coordindex, y=cg_rel]{\datay};
	\addplot[color=black,fill=colorx]	table[x expr=\coordindex, y=cg_rel]{\datayy};
    \end{groupplot}
    \node(title) at (3.7,2.8) {\large Normalized PCG counts, $n_{cg} / n_{cg}^{\text{first-order}}$};
	\end{tikzpicture}
	\caption{Normalized PCG iteration counts on all tested SuiteSparse SPD matrices 
	with at least 500,000 rows. For a given matrix and accuracy parameter $\varepsilon$, the number of iterations needed for convergence using the second-order schemes, is almost exactly halved compared to the first-order scheme.
}
	\label{fig:cg-suite-relative}
\end{figure}

\pgfplotstableread[col sep=comma]{logs/suite-sparse-es/prescale/hook_convergence_FO_0.01.csv}\datahookfo
\pgfplotstableread[col sep=comma]{logs/suite-sparse-es/prescale/hook_convergence_SO_0.01.csv}\datahookso
\pgfplotstableread[col sep=comma]{logs/suite-sparse-es/prescale/bone010_convergence_FO_0.01.csv}\databonefo
\pgfplotstableread[col sep=comma]{logs/suite-sparse-es/prescale/bone010_convergence_SO_0.01.csv}\databoneso
\pgfplotstableread[col sep=comma]{logs/suite-sparse-es/prescale/flan_convergence_FO_0.01.csv}\dataflanfo
\pgfplotstableread[col sep=comma]{logs/suite-sparse-es/prescale/flan_convergence_SO_0.01.csv}\dataflanso
\pgfplotstableread[col sep=comma]{logs/suite-sparse-es/prescale/audikw_convergence_FO_0.01.csv}\dataaudikwfo
\pgfplotstableread[col sep=comma]{logs/suite-sparse-es/prescale/audikw_convergence_SO_0.01.csv}\dataaudikwso

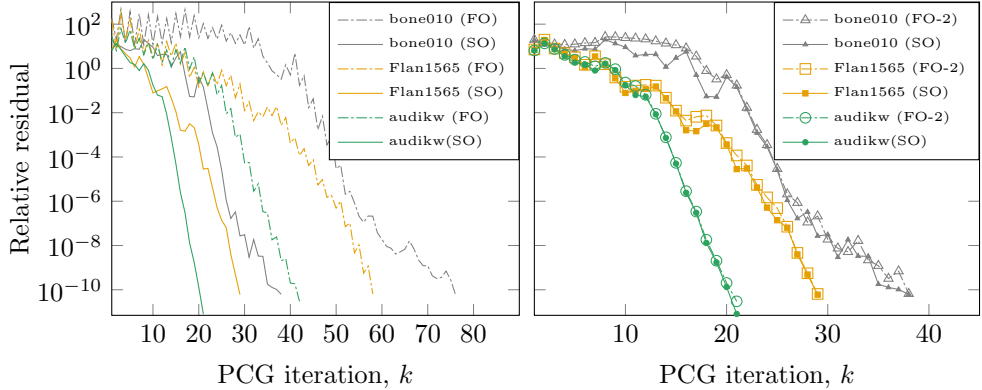
\begin{figure}[htbp]
    \centering
    \begin{tikzpicture}
        \begin{groupplot}[
            group style={
                group name=cg_converge,
                group size=2 by 1,
                xlabels at=edge bottom,
                xticklabels at=edge bottom,
                vertical sep=0.2cm,
                horizontal sep=0.2cm,
                ylabels at=edge left,
                yticklabels at=edge left,
            },
            ymode=log,ymin=7*10^(-12),ymax=1000.0,
            ytick={10^2,10^0,10^(-2),10^(-4),10^(-6),10^(-8),10^(-10)},
            yticklabels={$10^2$,$10^0$,$10^{-2}$,$10^{-4}$,$10^{-6}$,$10^{-8}$,$10^{-10}$},
            xmin=1,xmax=90,
            xtick={10,20,30,40,50,60,70,80},
            xticklabels={10,20,30,40,50,60,70,80},
        ]
        \nextgroupplot[width=7cm,height=5.75cm,xlabel={PCG iteration, $k$},
        	ylabel style={align=center},
        	ylabel={Relative residual},
        	legend columns = 1,
        	legend entries={bone010 (FO), bone010 (SO),Flan1565 (FO),Flan1565 (SO),audikw (FO), audikw(SO)},
        	legend style={font=\tiny},
        	legend cell align={left},
        	legend style={at={(1.0,1.0)}},
        	]
             \addplot[convfo1g_nomark] table[x=it, y=res]{\databonefo};
             \addplot[convso1g_nomark] table[x=it, y=res]{\databoneso};
             \addplot[convfo2b_nomark] table[x=it, y=res]{\dataflanfo};
             \addplot[convso2b_nomark] table[x=it, y=res]{\dataflanso};
             \addplot[convfo3_nomark] table[x=it, y=res]{\dataaudikwfo};
             \addplot[convso3_nomark] table[x=it, y=res]{\dataaudikwso};
       
         \nextgroupplot[width=7.5cm,height=5.75cm,xmax=45,xlabel={PCG iteration, $k$},
        	legend columns = 1,
        	legend entries={bone010 (FO-2), bone010 (SO),Flan1565 (FO-2),Flan1565 (SO),audikw (FO-2), audikw(SO)},
        	legend style={font=\tiny},
        	legend cell align={left},
        	legend style={at={(1.0,1.0)}},
        	]
             \addplot[convfo1g] table[x=it, y=res-2]{\databonefo};
             \addplot[convso1g] table[x=it, y=res]{\databoneso};
             \addplot[convfo2b] table[x=it, y=res-2]{\dataflanfo};
             \addplot[convso2b] table[x=it, y=res]{\dataflanso};
             \addplot[convfo3] table[x=it, y=res-2]{\dataaudikwfo};
             \addplot[convso3] table[x=it, y=res]{\dataaudikwso};
        \end{groupplot}
    \end{tikzpicture}
    \caption{Error decay in PCG preconditioned with spaND (accuracy $\varepsilon = 0.01$), for selected SuiteSparse matrices. On the left, we compare the convergence using the \fo{} (FO) and the \fullsec{} (SO).
    On the right, for each iteration $k$ we compare the residual obtained when using the \fullsec{} (SO) with the residual at the $2k$-th iteration when using the \fo{} (FO-2).
}
    \label{fig:cg-convergence}
\end{figure}

\begin{figure}[htbp]
	\begin{tikzpicture}
	\begin{groupplot}[
		group style={
        group name=suite-sparse-1,
        group size=1 by 2,
        xlabels at=edge bottom,
        xticklabels at=edge bottom,
        vertical sep=0cm,
        ylabels at=edge left,
        yticklabels at=edge left,
        vertical sep=0.2cm
    },
    xmin=-1,
    xmax=32,
    ymin=0,
    ymax=2,
    ybar=0.5pt,
    ytick={1.0,1.5,2.0},
    yticklabels={1.0,1.5,2.0},
    yticklabel style={rotate=-90},
    width=0.95\textwidth,
    xtick=data,
    xticklabels from table={\datay}{matrix},
    xticklabel style={tick label style={rotate=90}, font=\scriptsize},
    xtick={0,1,2,3,4,5,6,7,8,9,10,11,12,13,14,15,16,17,18,19,20,21,22,23,24,25,26,27,28,29,30,31},
    ]
    \nextgroupplot[
    	bar width=0.08cm,
    	height=3.75cm,
    	width=13.0cm,
    	legend entries={First-Order, Second-Order-Full, Second-Order-SF},
    	legend style={at={(1.0,1.2)},
    							font=\small},
    	legend columns=3,
    	ylabel style={align=left},
    	ylabel={$\varepsilon=0.1$},
    	ymajorgrids,
    ]
	\addplot[color=black,fill=none]	table[x expr=\coordindex, y expr=1.0]{\datax};
	\addplot[color=black,fill=colorz]	table[x expr=\coordindex, y=nnzfact_rel]{\datax};
	\addplot[color=black,fill=colorx]	table[x expr=\coordindex, y=nnzfact_rel]{\dataxx};
	    \nextgroupplot[
    	bar width=0.08cm,
    	height=3.75cm,
    	width=13.0cm,
    	ylabel style={align=left},
    	ylabel={$\varepsilon=0.05$},
    	ymajorgrids,
    ]
	\addplot[color=black,fill=none]	table[x expr=\coordindex, y expr=1.0]{\datay};
	\addplot[color=black,fill=colorz]	table[x expr=\coordindex, y=nnzfact_rel]{\datay};
	\addplot[color=black,fill=colorx]	table[x expr=\coordindex, y=nnzfact_rel]{\datayy};
    \end{groupplot}
    \node(title) at (3.7,2.9) {\large Normalized memory, $\mu/\mu^{\text{first-order}}$};
	\end{tikzpicture}
	\caption{Normalized memory requirements for the SuiteSparse matrices. The increase in memory requirements using the second-order schemes, never exceeds 100\%.
	When using the \superfine, it never exceeds 50\%.
}
\label{fig:density-suite-relative}
\end{figure}
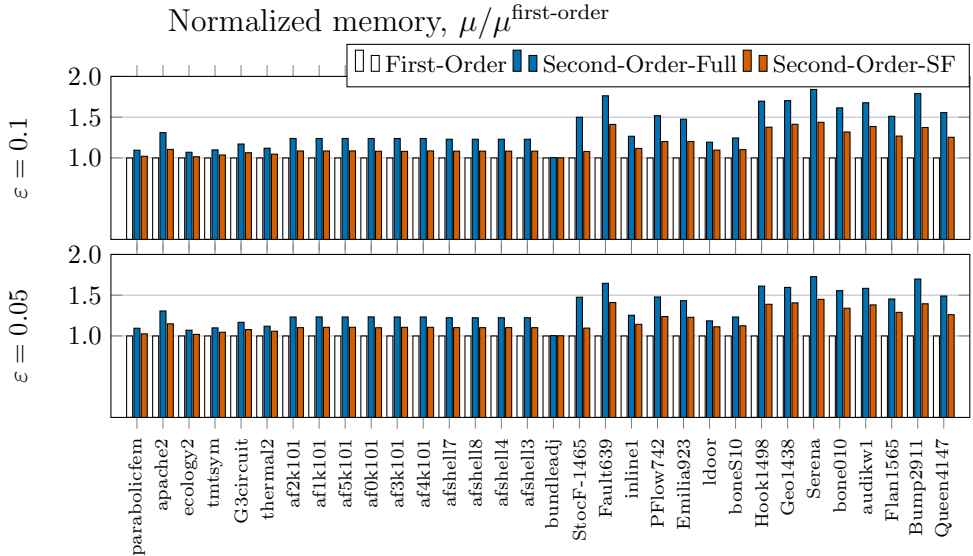

\begin{figure}[htbp]
	\begin{tikzpicture}
	\begin{groupplot}[
		group style={
        group name=suite-sparse-1,
        group size=1 by 2,
        xlabels at=edge bottom,
        xticklabels at=edge bottom,
        vertical sep=0cm,
        ylabels at=edge left,
        yticklabels at=edge left,
        vertical sep=0.2cm
    },
    xmin=-1,
    xmax=32,
    ymin=0,
    ymax=1.1,
    ybar=0.5pt,
    ytick={0.5,1.0},
    yticklabels={0.5,1.0},
    yticklabel style={rotate=-90},
    width=0.95\textwidth,
    xtick=data,
    xticklabels from table={\datay}{matrix},
    xticklabel style={tick label style={rotate=90}, font=\scriptsize},
    xtick={0,1,2,3,4,5,6,7,8,9,10,11,12,13,14,15,16,17,18,19,20,21,22,23,24,25,26,27,28,29,30,31},
    ]
    \nextgroupplot[
    	bar width=0.08cm,
    	height=4.0cm,
    	width=13.0cm,
    	legend entries={First-Order, Second-Order-Full, Second-Order-SF},
    	legend style={at={(1.0,1.45)},
    							font=\small},
    	ylabel style={align=left},
    	ylabel={$\varepsilon=0.1$},
    	ymajorgrids,
    ]
	\addplot[color=black,fill=none]	table[x expr=\coordindex, y expr=1.0]{\datax};
	\addplot[color=black,fill=colorz]	table[x expr=\coordindex, y=tsolve_rel]{\datax};
	\addplot[color=black,fill=colorx]	table[x expr=\coordindex, y=tsolve_rel]{\dataxx};
	    \nextgroupplot[
    	bar width=0.08cm,
    	height=4.0cm,
    	width=13.0cm,
    	ylabel style={align=left},
    	ylabel={$\varepsilon=0.05$},
    	ymajorgrids,
    ]
	\addplot[color=black,fill=none]	table[x expr=\coordindex, y expr=1.0]{\datay};
	\addplot[color=black,fill=colorz]	table[x expr=\coordindex, y=tsolve_rel]{\datay};
	\addplot[color=black,fill=colorx]	table[x expr=\coordindex, y=tsolve_rel]{\datayy};
    \end{groupplot}
    \node(title) at (3.7,2.8) {\large Normalized solve times, $T_s/T_s^{\text{first-order}}$};
	\end{tikzpicture}
	\caption{Normalized solve times for the SuiteSparse matrices. For a given matrix and accuracy parameter $\varepsilon$, the time needed for convergence of PCG is significantly reduced when using the second-order schemes.}
\label{fig:solve-suite-relative}
\end{figure}
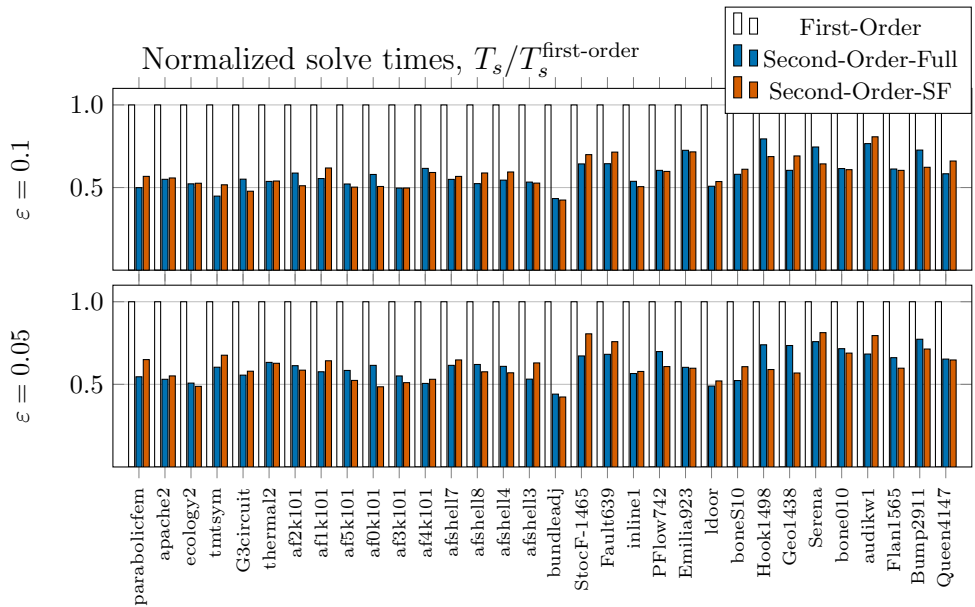

\makeatletter
\newcommand\resetstackedplots{
\makeatletter
\pgfplots@stacked@isfirstplottrue
\makeatother
\addplot [forget plot,draw=none] coordinates{(0,0) (1,0) (2,0) (3,0) (4,0) (5,0) (6,0) (7,0) (8,0) (9,0) (10,0) (11,0) (12,0) (13,0) (14,0) (15,0) (16,0) (17,0) (18,0) (19,0) (21,0) (22,0) (23,0) (24,0) (25,0) (25,0) (27,0) (28,0) (29,0) (30,0) (31,0) (32,0)};
}
\makeatother

\pgfplotstableread[col sep=comma]{logs/suitesparse-3-table-FO-prescale-rel.csv}\dataxrel
\pgfplotstableread[col sep=comma]{logs/suitesparse-3-table-SO-prescale-rel.csv}\dataxxrel
\pgfplotstableread[col sep=comma]{logs/suitesparse-3-table-SO-sf-prescale-rel.csv}\dataxxxrel
\begin{figure}[htbp]
	\begin{tikzpicture}
	\begin{groupplot}[
		group style={
        group name=suite-sparse-1,
        group size=1 by 1,
        xlabels at=edge bottom,
        xticklabels at=edge bottom,
        vertical sep=0cm,
        ylabels at=edge left,
        yticklabels at=edge left,
        vertical sep=0.2cm
    },
    xmin=-1,
    xmax=32,
    ymin=0,
    ymax=1,
    ybar stacked,
    ytick={0.25,0.5,0.75,1.0},
    yticklabels={0.25,0.5,0.75,1.0},
    yticklabel style={rotate=-90},
    width=0.95\textwidth,
    xtick=data,
    xticklabels from table={\datay}{matrix},
    xticklabel style={tick label style={rotate=90}, font=\scriptsize},
    xtick={0,1,2,3,4,5,6,7,8,9,10,11,12,13,14,15,16,17,18,19,20,21,22,23,24,25,26,27,28,29,30,31,32},
    ]
    \nextgroupplot[
    	bar width=0.08cm,
    	height=5.0cm,
    	width=13.5cm,
    	legend entries={Factorization (FO), PCG Solve (FO), Factorization (SO-Full), PCG Solve (SO-Full), Factorization (SO-SF), PCG Solve (SO-SF)},
    	legend columns=2,
    	transpose legend,
    	legend style={at={(1.0,1.3)},
    							font=\footnotesize},
    	ymajorgrids,
    ]
	\addplot[color=black,fill=gray!50]	table[x expr=\coordindex-0.3, y=tf_rel]{\dataxrel};
	\addplot[color=black,fill=none]	table[x expr=\coordindex-0.3, y=ts_rel]{\dataxrel};
	\resetstackedplots
	
	\addplot[color=black,fill=colorz]	table[x expr=\coordindex, y=tf_rel]{\dataxxrel};
	\addplot[color=black,fill=colorz!50]	table[x expr=\coordindex, y=ts_rel]{\dataxxrel};
	\resetstackedplots

	\addplot[color=black,fill=colorx]	table[x expr=\coordindex+0.3, y=tf_rel]{\dataxxxrel};
	\addplot[color=black,fill=colorx!50]	table[x expr=\coordindex+0.3, y=ts_rel]{\dataxxxrel};
		
  \end{groupplot}
     \node(title) at (5.5 ,4.75) {\large Normalized total runtimes (optimal $\varepsilon$ for each problem and scheme)};
	\end{tikzpicture}
	\caption{Normalized total runtimes ($T_t/T_t^{\text{first-order}} = T_f/T_t^{\text{first-order}} + T_s/T_t^{\text{first-order}}$) for the SuiteSprase matrices. For each problem, and each tested scheme, the optimal corresponding $\varepsilon$ was chosen, from the four tested ones ($\varepsilon=0.2,0.1,0.05,0.01$). Second-order schemes improve the total runtimes in all tested cases.}
	\label{fig:cg-suite-relative-2}
\end{figure}
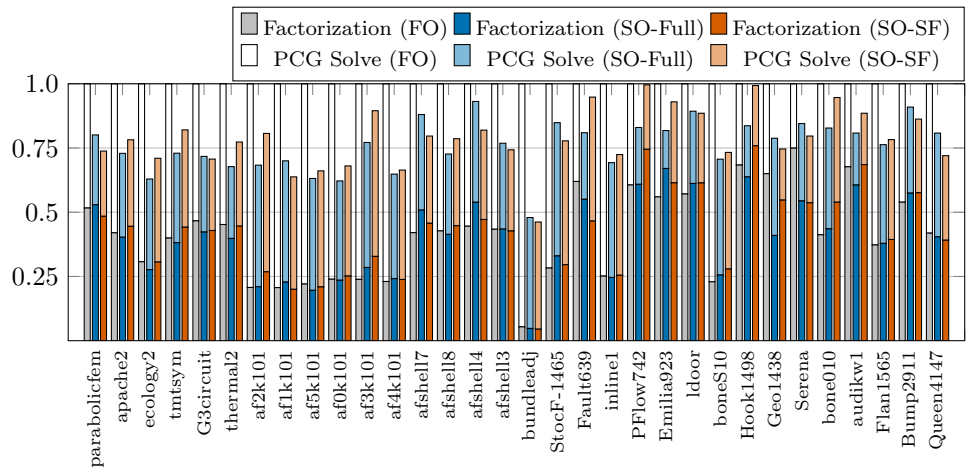


\subsubsection{Improved total timings}
Since for a given $\varepsilon$ parameter the factorization phase of spaND involves little additional computations when using the second-order schemes, the total runtimes should also be reduced.
In \Cref{fig:cg-suite-relative-2}, for each tested matrix, and each scheme, we show the best runtime from amongst the four tested accuracy parameters $\varepsilon$, split into the factorization and solve (PCG iteration) phases (full data can be found in \Cref{sec:appendix}).
We observe improvements in the optimal total runtime in all tested cases when using the new second-order schemes.

Notice that a significantly shorter solve time $T_s$ for the same accuracy parameter $\varepsilon$, may mean that---from the standpoint of optimizing the total runtime---the factorization is too accurate.
Therefore one can expect that the optimal total runtime when using the second-order schemes will be obtained for a larger accuracy parameter $\varepsilon$, than when using the \fo.
We observe this for a number of test cases (for example $\texttt{Geo1436}$ or $\texttt{Serena}$, see also \Cref{sec:appendix}).

\subsection{Differences between the two second-order scheme variants}
For the same accuracy parameter $\varepsilon$, on almost all tested problems, the \superfine{} and the \fullsec{} resulted in nearly the same PCG iteration counts. The \superfine{} preconditioner has lower memory requirements than the \fullsec{} preconditioner, and is cheaper to apply.
This may translate to savings in the solution phase.
The rank-revealing QR has to be computed down to $\varepsilon^2$ accuracy, however, 
which may make the factorization phase more expensive.
Both second-order scheme variants performed competitively in our test cases.

%% file: conclusions.tex
\section{Conclusions}
\label{sec:conclusions}

We introduced a second-order accurate approach to sparsifying the numerically low-rank blocks in the approximate hierarchical factorizations of sparse symmetric positive definite matrices.
Similar to the standard first-order approach, we apply orthogonal matrices defined by the rank-revealing decomposition of the given off-diagonal block, so that interactions of many variables become small, and are subsequently dropped.
However, the new approach also includes additional terms that approximately eliminate these variables.
As a result, the norm of the overall error depends quadratically, as opposed to linearly, on the norm of the error in the low-rank approximation of the given block.

Numerical analysis of the resulting two-level preconditioners, as well as numerical experiments, show clear improvements provided by the new method.
In particular, the analysis suggests that the number of Conjugate Gradient iterations should be halved  for any given accuracy parameter.
Consistent with this, when incorporated into the spaND algorithm \cite{cambier2020algebraic}, for any given accuracy, the new approach results in a reduction of iteration counts by almost exactly half, on a wide range of SPD problems.
The new approach involves little additional computations in the factorization phase, and improves the total runtimes of spaND.

Beside spaND, other solvers based on hierarchical low-rank structures can benefit from our results when applied to sparse matrices.
In particular, we considered only factorizations sparsifying all off-diagonal blocks but the second-order scheme can be similarly defined for algorithms distinguishing neighboring and well-separated interactions, such as \cite{sushnikova2018compress,pouransari2017fast,minden2017recursive,chen2018distributed}.
Also, the sparsification approach that improves accuracy on the chosen near-kernel subspace, as in \cite{klockiewicz2019sparse}, can be applied basically without modifications.
Lastly, the new approach can be expected to work optimally for a lower accuracy parameter, as observed on some test problems.
This should improve the parallel properties of hierarchical solvers because a larger portion of computations is then performed on small blocks in the initial levels of the algorithm.

%% file: acknowledgements.tex
\section*{Acknowledgments}
Bazyli Klockiewicz would like to thank the Stanford University Petroleum Research Institute's Reservoir Simulation Industrial Affiliates Program (SUPRI-B) for its financial support. L\'eopold Cambier was supported by a Fellowship from Total S.A.
 Some of the computing was performed on the Sherlock cluster at Stanford University. We would like to thank Stanford University and the Stanford Research Computing Center for providing computational resources and support.
 We would also like to thank Jordi Feliu-Fab\`{a} for useful conversations about the numerical analysis of the algorithms, and Abeynaya Gnanasekaran for providing the implementation of the early-stopping column-pivoted QR.

%% file: appendix.tex
\newpage
\appendix
\section{Proof of \Cref{thm:residuals}}
\label{sec:proof}
\newline
\begin{proof}
Denote $F = I - C^\top C.$ We have
$
A_1
\begin{pmatrix}
	b_1 \\ b_2
\end{pmatrix} =
\begin{pmatrix}
	Cb_2 + b_1 \\ b_2
\end{pmatrix}
$. We also have \[
\begin{pmatrix}
	-Cb_2 + b_1 \\ b_2
\end{pmatrix}
=
-
\begin{pmatrix}
	Cb_2 + b_1 \\ b_2
\end{pmatrix}
+ 2
\begin{pmatrix}
	b_1 \\ b_2
\end{pmatrix}
\] 
and
\begin{equation}
\label{eq:a1-krylov}
A_1 
\begin{pmatrix}
	-Cb_2 + b_1 \\ b_2
\end{pmatrix}
=
\begin{pmatrix}
	b_1 \\ Fb_2
\end{pmatrix}
\end{equation}
By iterating this (with $Fb_2$ in place of $b_2$), we have that one particular sequence of vectors defining the Krylov subspaces for $(A_1,b)$, is given by:
\begin{gather}
\mathcal{K}(A_1,b) \\
= \nonumber
\Big\{
\begin{pmatrix}
b_1 \\ b_2
\end{pmatrix},
\begin{pmatrix}
-Cb_2 + b_1 \\ b_2
\end{pmatrix}, 
\begin{pmatrix}
b_1 \\ F b_2
\end{pmatrix}, 
\begin{pmatrix}
-CFb_2 + b_1 \\ F b_2
\end{pmatrix},
\begin{pmatrix}
b_1 \\ F^2 b_2
\end{pmatrix},
\begin{pmatrix}
-CF^2b_2 + b_1 \\ F^2b_2
\end{pmatrix},
\dots \Big\}
\end{gather}
An analogous sequence for $(A_2,b)$ is given by
\[
\mathcal{K}(A_2,b) =
\Big\{
\begin{pmatrix}
b_1 \\ b_2
\end{pmatrix},
\begin{pmatrix}
b_1 \\ Fb_2
\end{pmatrix},
\begin{pmatrix}
b_1 \\ F^2b_2
\end{pmatrix},
\dots
\Big\}
\]
There exist coefficients $\alpha_i$, for $0 \leq i \leq k-1,$ such that the Conjugate Gradient solution to $A_2 x = b$ at the $k$-th iteration, is given by:
\[
x^{(2)}_k =
\sum_{i=0}^{k-1} \alpha_i 
\begin{pmatrix}
b_1 \\ F^i b_2
\end{pmatrix}
 \]
The residual then equals
 \[
r^{(2)}_k = b - A_2 \sum_{i=0}^{k-1} \alpha_i 
\begin{pmatrix}
b_1 \\ F^i b_2
\end{pmatrix}
=
b - \sum_{i=0}^{k-1} \alpha_i 
\begin{pmatrix}
b_1 \\ F^{i+1} b_2
\end{pmatrix}
\]
From the optimality conditions of the Conjugate Gradient solution, we have
\begin{equation} \label{eq:perp-1}
r^{(2)}_k =
b - \sum_{i=0}^{k-1} \alpha_i 
\begin{pmatrix}
b_1 \\ F^{i+1} b_2
\end{pmatrix}
\perp 
\langle
\Big\{ 
\begin{pmatrix}
b_1 \\ F^i b_2
\end{pmatrix}
\Big\}_{i=0}^{k-1}
\rangle
\end{equation}
Define:
\[
s:=
b - 
A_1
\sum_{i=0}^{k-1}
\alpha_i
\begin{pmatrix}
-CF^i b_2 + b_1 \\ F^i b_2
\end{pmatrix}
\]
From \ref{eq:a1-krylov} and \ref{eq:perp-1}, we have:
\begin{equation}\label{eq:perp-2}
s =
b - \sum_{i=0}^{k-1} \alpha_i 
\begin{pmatrix}
b_1 \\ F^{i+1} b_2
\end{pmatrix}
=
r_k^{(2)}
\perp 
\langle
\Big\{ 
\begin{pmatrix}
b_1 \\ F^i b_2
\end{pmatrix}
\Big\}_{i=0}^{k-1}
\rangle
\end{equation}
We now show that also
\[
s 
\perp 
\langle
\Big\{ 
\begin{pmatrix}
-CF^i b_2 +b_1 \\ F^i b_2
\end{pmatrix}
\Big\}_{i=0}^{k-1}
\rangle
\]
From \ref{eq:perp-2} it suffices to show that
\[
s =
b -
\sum_{i=0}^{k-1} \alpha_i 
\begin{pmatrix}
b_1 \\ F^{i+1} b_2
\end{pmatrix}
\perp 
\langle
\Big\{ 
\begin{pmatrix}
-CF^i b_2 \\ 0
\end{pmatrix}
\Big\}_{i=0}^{k-1}
\rangle
\]
However, this holds, because for any $0 \leq i \leq k-1$
\[
s^\top \begin{pmatrix}
-CF^i b_2 \\ 0
\end{pmatrix} = \Big(\sum_{j=0}^{k-1}{\alpha_j} - 1\Big)b_1^\top CF^ib_2 = 0
\]
which follows from the assumption that $b_1^\top C = 0$.
This means that $s$ is orthogonal to the $2k$-th Krylov subspace for $A_1$, and therefore, again from the optimality conditions
\[
x^{(1)}_{2k} := \sum_{i=0}^{k-1} \alpha_i 
\begin{pmatrix}
-CF^i b_2 + b_1 \\ F^i b_2
\end{pmatrix}
\]
is the solution to $A_1 x = b$ produced at the $2k$-th iteration, and
\[
r^{(1)}_{2k} = s = r_k^{(2)}
\]
\end{proof}

\section{Tables with runs for optimal $\varepsilon$}
\label{sec:appendix}
For each tested matrix from the SuiteSparse collection, and each sparsification scheme, we show the run with optimal $\varepsilon$ (from among the four tested, $\varepsilon = 0.2, 0.1, 0.05, 0.01$) in terms of the total runtime.
Some of the matrices (almost identical to already shown), were omitted.
\bgroup
\setlength{\tabcolsep}{0.35em}
\begin{table}[htbp]
\renewcommand{\arraystretch}{1.1}
\small
\centering {
	\begin{tabular}{llll | rrrrrrr >{\bfseries}r}
	 Matrix &  $n$ & $nnz(A)$ &Order &  $\varepsilon$ & $\mu$ & $n_{cg}$ &$T_f[s]$ & $T_s[s]$ & $T_t[s]$ & $\Delta T_{t}$ \\
	 \hline
parabolicfem&$5.3\cdot 10^5$&$3.7\cdot 10^6$&First&0.01&8.9&10&2.9&2.7&5.6&--\\
&$5.3\cdot 10^5$&$3.7\cdot 10^6$&Sec-Full&0.01&9.8&5&3.0&1.5&4.5&\textbf{-19.9\%}\\
&$5.3\cdot 10^5$&$3.7\cdot 10^6$&Sec-SF&0.01&9.2&5&2.7&1.4&4.1&\textbf{-26.2\%}\\
	 \hline
apache2&$7.2\cdot 10^5$&$4.8\cdot 10^6$&First&0.01&17.0&24&11.9&16.5&28.4&--\\
&$7.2\cdot 10^5$&$4.8\cdot 10^6$&Sec-Full&0.01&21.9&12&11.5&9.3&20.7&\textbf{-27.1\%}\\
&$7.2\cdot 10^5$&$4.8\cdot 10^6$&Sec-SF&0.01&20.2&12&12.7&9.6&22.2&\textbf{-21.8\%}\\
	 \hline
ecology2&$1.0\cdot 10^6$&$5.0\cdot 10^6$&First&0.01&10.3&17&2.5&5.7&8.2&--\\
&$1.0\cdot 10^6$&$5.0\cdot 10^6$&Sec-Full&0.01&11.0&8&2.3&2.9&5.1&\textbf{-37.1\%}\\
&$1.0\cdot 10^6$&$5.0\cdot 10^6$&Sec-SF&0.01&10.6&9&2.5&3.3&5.8&\textbf{-29.0\%}\\
	 \hline
tmtsym&$7.3\cdot 10^5$&$5.1\cdot 10^6$&First&0.01&6.9&20&5.3&7.9&13.2&--\\
&$7.3\cdot 10^5$&$5.1\cdot 10^6$&Sec-Full&0.01&7.6&10&5.0&4.6&9.6&\textbf{-27.0\%}\\
&$7.3\cdot 10^5$&$5.1\cdot 10^6$&Sec-SF&0.01&7.3&10&5.8&5.0&10.8&\textbf{-17.9\%}\\
	 \hline
G3circuit&$1.6\cdot 10^6$&$7.7\cdot 10^6$&First&0.01&12.9&14&12.2&14.0&26.2&--\\
&$1.6\cdot 10^6$&$7.7\cdot 10^6$&Sec-Full&0.01&14.9&7&11.1&7.7&18.8&\textbf{-28.3\%}\\
&$1.6\cdot 10^6$&$7.7\cdot 10^6$&Sec-SF&0.01&14.1&7&11.3&7.3&18.6&\textbf{-29.3\%}\\
	 \hline
thermal2&$1.2\cdot 10^6$&$8.6\cdot 10^6$&First&0.01&6.3&16&10.2&12.4&22.6&--\\
&$1.2\cdot 10^6$&$8.6\cdot 10^6$&Sec-Full&0.01&7.0&8&9.0&6.3&15.3&\textbf{-32.3\%}\\
&$1.2\cdot 10^6$&$8.6\cdot 10^6$&Sec-SF&0.01&6.7&8&10.1&7.4&17.5&\textbf{-22.7\%}\\
\end{tabular}
	}
\caption{(Part 1/3) SuiteSparse SPD matrices: best runs in terms of the total runtime.
 }
\label{tab:suite-best-1}
\end{table}
\egroup
\bgroup
\setlength{\tabcolsep}{0.35em}
\begin{table}[htbp]
\renewcommand{\arraystretch}{1.1}
\small
\centering {
	\begin{tabular}{llll | rrrrrrr >{\bfseries}r}
	 Matrix &  $n$ & $nnz(A)$ &Order &  $\varepsilon$ & $\mu$ & $n_{cg}$ &$T_f[s]$ & $T_s[s]$ & $T_t[s]$ & $\Delta T_{t}$ \\
	 \hline
afshell8&$5.0\cdot 10^5$&$1.8\cdot 10^7$&First&0.01&3.9&16&3.9&5.0&8.9&--\\
&$5.0\cdot 10^5$&$1.8\cdot 10^7$&Sec-Full&0.01&4.7&8&3.9&3.0&6.8&\textbf{-23.1\%}\\
&$5.0\cdot 10^5$&$1.8\cdot 10^7$&Sec-SF&0.01&4.4&8&3.8&2.8&6.6&\textbf{-25.7\%}\\
	 \hline
bundleadj&$5.1\cdot 10^5$&$2.0\cdot 10^7$&First&0.01&2.3&86&0.8&13.4&14.1&--\\
&$5.1\cdot 10^5$&$2.0\cdot 10^7$&Sec-Full&0.01&2.3&42&0.7&6.1&6.8&\textbf{-52.0\%}\\
&$5.1\cdot 10^5$&$2.0\cdot 10^7$&Sec-SF&0.01&2.3&45&0.6&5.9&6.5&\textbf{-53.8\%}\\
	 \hline
StocF-1465&$1.5\cdot 10^6$&$2.1\cdot 10^7$&First&0.01&20.7&122&137.3&348.3&485.6&--\\
&$1.5\cdot 10^6$&$2.1\cdot 10^7$&Sec-Full&0.01&29.9&58&160.2&251.8&412.1&\textbf{-15.1\%}\\
&$1.5\cdot 10^6$&$2.1\cdot 10^7$&Sec-SF&0.01&24.1&61&143.7&234.0&377.7&\textbf{-22.2\%}\\
	 \hline
Fault639&$6.4\cdot 10^5$&$2.9\cdot 10^7$&First&0.05&14.1&44&106.9&65.6&172.5&--\\
&$6.4\cdot 10^5$&$2.9\cdot 10^7$&Sec-Full&0.05&23.1&21&94.9&44.7&139.6&\textbf{-19.1\%}\\
&$6.4\cdot 10^5$&$2.9\cdot 10^7$&Sec-SF&0.1&16.3&46&80.3&83.2&163.5&\textbf{-5.2\%}\\
\hline
af5k101&$5.0\cdot 10^5$&$1.8\cdot 10^7$&First&0.01&3.9&40&3.8&12.8&16.7&--\\
&$5.0\cdot 10^5$&$1.8\cdot 10^7$&Sec-Full&0.01&4.8&19&4.0&6.8&10.8&\textbf{-35.2\%}\\
&$5.0\cdot 10^5$&$1.8\cdot 10^7$&Sec-SF&0.01&4.4&20&4.0&7.1&11.1&\textbf{-33.6\%}\\
\hline
inline1&$5.0\cdot 10^5$&$3.7\cdot 10^7$&First&0.01&3.7&72&12.8&38.0&51.2&--\\
&$5.0\cdot 10^5$&$3.7\cdot 10^7$&Sec-Full&0.01&4.5&33&12.5&22.7&35.2&\textbf{-30.7\%}\\
&$5.0\cdot 10^5$&$3.7\cdot 10^7$&Sec-SF&0.01&4.3&33&12.9&23.9&36.8&\textbf{-27.5\%}\\
	 \hline
PFlow742&$7.4 \cdot 10^5$&$3.7 \cdot 10^7$&First&0.01&6.9&22&39.8&25.9&65.7&--\\
&$7.4 \cdot 10^5$&$3.7 \cdot 10^7$&Sec-Full&0.01&9.7&10&40.0&14.5&54.5&\textbf{-17.0\%}\\
&$7.4 \cdot 10^5$&$3.7 \cdot 10^7$&Sec-SF&0.01&8.9&10&48.9&16.5&65.4&\textbf{-0.4\%}\\
%
\hline
Emilia923&$9.2 \cdot 10^5$&$4.1 \cdot 10^7$&First&0.1&17.6&80&199.1&156.6&355.7&--\\
&$9.2 \cdot 10^5$&$4.1 \cdot 10^7$&Sec-Full&0.05&27.7&17&238.4&52.5&291.0&\textbf{-18.2\%}\\
&$9.2 \cdot 10^5$&$4.1 \cdot 10^7$&Sec-SF&0.1&21.1&41&218.6&112.1&330.7&\textbf{-7.0\%}\\
	 \hline
ldoor&$9.5 \cdot 10^5$&$4.7 \cdot 10^7$&First&0.01&2.8&9&7.4&5.6&13.0&--\\
&$9.5 \cdot 10^5$&$4.7 \cdot 10^7$&Sec-Full&0.01&3.3&5&7.9&3.7&11.6&\textbf{-10.7\%}\\
&$9.5 \cdot 10^5$&$4.7 \cdot 10^7$&Sec-SF&0.01&3.1&5&8.0&3.5&11.5&\textbf{-11.5\%}\\
	 \hline
boneS10&$9.1 \cdot 10^5$&$5.5 \cdot 10^7$&First&0.01&4.0&64&16.5&55.3&71.8&--\\
&$9.1 \cdot 10^5$&$5.5 \cdot 10^7$&Sec-Full&0.01&4.8&31&18.3&32.4&50.7&\textbf{-29.3\%}\\
&$9.1 \cdot 10^5$&$5.5 \cdot 10^7$&Sec-SF&0.01&4.5&31&20.0&32.6&52.6&\textbf{-26.7\%}\\
	 \hline
Hook1498&$1.5 \cdot 10^6$&$6.1 \cdot 10^7$&First&0.01&13.7&34&209.6&97.0&306.5&--\\
&$1.5 \cdot 10^6$&$6.1 \cdot 10^7$&Sec-Full&0.01&20.1&17&195.4&61.1&256.4&\textbf{-16.3\%}\\
&$1.5 \cdot 10^6$&$6.1 \cdot 10^7$&Sec-SF&0.01&18.4&17&232.5&71.9&304.4&\textbf{-0.7\%}\\
	 \hline
Geo1438&$1.4 \cdot 10^6$&$6.3 \cdot 10^7$&First&0.05&15.0&39&263.6&141.9&405.6&--\\
&$1.4 \cdot 10^6$&$6.3 \cdot 10^7$&Sec-Full&0.1&20.9&37&166.1&153.5&319.6&\textbf{-21.2\%}\\
&$1.4 \cdot 10^6$&$6.3 \cdot 10^7$&Sec-SF&0.05&21.1&19&221.9&80.7&302.6&\textbf{-25.4\%}\\
 \hline
 Serena&$1.4 \cdot 10^6$&$6.5 \cdot 10^7$&First&0.05&13.7&31&292.1&97.7&389.8&--\\
&$1.4 \cdot 10^6$&$6.5 \cdot 10^7$&Sec-Full&0.1&20.7&25&212.1&117.4&329.5&\textbf{-15.5\%}\\
&$1.4 \cdot 10^6$&$6.5 \cdot 10^7$&Sec-SF&0.1&16.1&26&209.2&101.3&310.5&\textbf{-20.3\%}\\
	 \hline
bone010&$9.9 \cdot 10^5$&$7.2 \cdot 10^7$&First&0.01&8.5&76&104.6&149.0&253.6&--\\
&$9.9 \cdot 10^5$&$7.2 \cdot 10^7$&Sec-Full&0.01&12.2&38&110.4&99.4&209.9&\textbf{-17.2\%}\\
&$9.9 \cdot 10^5$&$7.2 \cdot 10^7$&Sec-SF&0.01&11.1&36&136.7&103.3&240.0&\textbf{-5.4\%}\\
\end{tabular}
	}
\caption{(Part 2/3)  SuiteSparse SPD matrices: best runs in terms of the total runtime.
 }
\label{tab:suite-best-3}
\end{table}
\egroup

\bgroup
\setlength{\tabcolsep}{0.35em}
\begin{table}[htbp]
\renewcommand{\arraystretch}{1.1}
\small
\centering {
	\begin{tabular}{llll | rrrrrrr >{\bfseries}r}
	 Matrix &  $n$ & $nnz(A)$ &Order &  $\varepsilon$ & $\mu$ & $n_{cg}$ &$T_f[s]$ & $T_s[s]$ & $T_t[s]$ & $\Delta T_{t}$ \\
	 \hline
audikw1&$9.4 \cdot 10^5$&$7.8 \cdot 10^7$&First&0.01&9.6&42&219.9&110.6&330.5&--\\
&$9.4 \cdot 10^5$&$7.8 \cdot 10^7$&Sec-Full&0.01&13.8&21&195.5&65.3&260.8&\textbf{-21.1\%}\\
&$9.4 \cdot 10^5$&$7.8 \cdot 10^7$&Sec-SF&0.01&12.7&21&221.2&64.5&285.7&\textbf{-13.5\%}\\
	 \hline
Flan1565&$1.6 \cdot 10^6$&$1.2 \cdot 10^8$&First&0.01&7.0&58&107.2&180.4&287.6&--\\
&$1.6 \cdot 10^6$&$1.2 \cdot 10^8$&Sec-Full&0.01&9.4&29&108.9&110.5&219.4&\textbf{-23.7\%}\\
&$1.6 \cdot 10^6$&$1.2 \cdot 10^8$&Sec-SF&0.01&8.8&29&113.3&111.8&225.1&\textbf{-21.7\%}\\
	 \hline
Bump2911&$2.9 \cdot 10^6$&$1.3 \cdot 10^8$&First&0.1&15.4&89&797.3&681.5&1478.8&--\\
&$2.9 \cdot 10^6$&$1.3 \cdot 10^8$&Sec-Full&0.1&27.5&39&849.3&495.3&1344.6&\textbf{-9.1\%}\\
&$2.9 \cdot 10^6$&$1.3 \cdot 10^8$&Sec-SF&0.1&21.1&42&851.2&424.3&1275.5&\textbf{-13.7\%}\\
	 \hline
Queen4147&$4.1 \cdot 10^6$&$3.3 \cdot 10^8$&First&0.2&9.5&110&842.0&1169.5&2011.5&--\\
&$4.1 \cdot 10^6$&$3.3 \cdot 10^8$&Sec-Full&0.2&15.7&49&812.8&812.6&1625.4&\textbf{-19.2\%}\\
&$4.1 \cdot 10^6$&$3.3 \cdot 10^8$&Sec-SF&0.2&11.6&57&786.6&661.9&1448.5&\textbf{-28.0\%}\\
\end{tabular}
	}
\caption{(Part 3/3)  SuiteSparse SPD matrices: best runs in terms of the total runtime.
 }
\label{tab:suite-best-4}
\end{table}
\egroup